\documentclass[reqno,11pt]{article}
\usepackage{amsmath, amsthm, amscd,amsfonts, amssymb,hyperref,floatrow}
\usepackage{graphicx, color,dsfont,xcolor,epsfig}

\numberwithin{equation}{section}

\newtheorem{theorem}{Theorem}[section]
\newtheorem{lemma}[theorem]{Lemma}
\newtheorem{proposition}[theorem]{Proposition}
\newtheorem{corollary}[theorem]{Corollary}
\newtheorem{rem}[theorem]{Remark}
\newtheorem{definition}[theorem]{Definition}

\newtheorem{assumption}[theorem]{Assumption}

\def\build#1_#2^#3{\mathrel{\mathop{\kern 0pt#1}\limits_{#2}^{#3}}}

 \hypersetup{
    linktoc=page,
    linkcolor=blue,          
    citecolor=red,        
    filecolor=blue,      
    urlcolor=cyan,
    colorlinks=true           
}

\renewcommand{\tilde}{\widetilde}          
\DeclareMathSymbol{\leqslant}{\mathalpha}{AMSa}{"36} 
\DeclareMathSymbol{\geqslant}{\mathalpha}{AMSa}{"3E} 
\DeclareMathSymbol{\eset}{\mathalpha}{AMSb}{"3F}     
\renewcommand{\leq}{\;\leqslant\;}                   
\renewcommand{\geq}{\;\geqslant\;}                   
\DeclareMathOperator{\sign}{sign}

\newcommand{\R}{\mathbb{R}}

\def \Cov{\mathrm{Cov}}

\def \E{ \mathbb E  }

\definecolor{remi}{rgb}{0,0,0}

\let\qed=\QED

\def\md{\mid}
\def \eps {\epsilon}
\def\Bb#1#2{{\def\md{\bigm| }#1\bigl[#2\bigr]}}

\def\Eb{\Bb\E}

\def\<#1{\langle #1\rangle}

\def\nn{\nonumber}
\def\bi{\begin{itemize}}  
\def\ei{\end{itemize}}
\def\bnum{\begin{enumerate}} 
\def\enum{\end{enumerate}}
\def\ni{\noindent}

\let\sc\scshape

\def\du{\delta_\ell u} 
\def\Pl{\Phi_\ell}

\title{On a skewed and multifractal uni-dimensional random field, as a probabilistic representation of Kolmogorov's views on turbulence}
\author{Laurent Chevillard\footnote{Univ Lyon, Ens de Lyon, Univ Claude Bernard, CNRS, Laboratoire de Physique, 46 all\'ee d'Italie F-69342 Lyon, France.},\, 
Christophe  Garban\footnote{Universit\'e de Lyon, Institut Camille Jordan, 43 blvd. du 11 novembre 1918, F-69622 Villeurbanne cedex, France. },\, 
R\'emi Rhodes\footnote{Universit\'e Paris Est-Marne la Vall\'ee, LAMA and CNRS UMR 8050, France.},\,
Vincent Vargas\footnote{ENS Ulm, DMA, 45 rue d'Ulm,  75005 Paris, France.}}

\begin{document}
\maketitle

\begin{abstract}
We construct, for the first time to our knowledge, a one-dimensional stochastic field $\{u(x)\}_{x\in \R}$ which satisfies the following axioms which are at the core of the phenomenology of turbulence mainly due to Kolmogorov:
\bnum
\item[(i)] Homogeneity and isotropy: $u(x) \overset{\mathrm{law}}= -u(x) \overset{\mathrm{law}}=u(0)$
\item[(ii)] Negative skewness (i.e. the $4/5^{\mbox{\tiny th}}$-law): \\
$\Eb{(u(x+\ell)-u(x))^3} \sim_{\ell \to 0+} - C \, \ell\,,$ \, for some constant $C>0$  
\item[(iii)] Intermittency: \\
$\Eb{|u(x+\ell)-u(x) |^q} \asymp_{\ell \to 0} |\ell|^{\xi_q}\,,$ \, for some non-linear spectrum $q\mapsto \xi_q$
\enum
Since then, it has been a challenging problem to combine axiom (ii) with axiom (iii) (especially for Hurst indexes of interest in turbulence, namely $H<1/2$). In order to achieve simultaneously both axioms, we disturb with two ingredients a underlying fractional Gaussian field of parameter $H\approx \frac 1 3 $. The first ingredient is an independent Gaussian multiplicative chaos (GMC) of parameter $\gamma$ that mimics the intermittent, i.e. multifractal, nature of the fluctuations. The second one is a field that correlates in an intricate way the fractional component and the GMC without additional parameters, a necessary inter-dependence in order to reproduce the asymmetrical, i.e. skewed, nature of the probability laws at small scales.
\end{abstract}

\tableofcontents

\vspace{1cm}
\footnotesize

\normalsize

\section{Introduction}

\subsection{General words concerning the present approach}

The quest for a rigorous probabilistic model of the velocity field in a 3d turbulent flow is a longstanding problem which goes back to the seminal work of Kolmogorov  \cite{Kol41}. The purpose of this work is to propose a new and tractable model in this direction. For the sake of simplicity, we will restrict to the simplified framework of 1d random fields. However, we believe our model can be generalized to the realistic 3d case: this generalization will be considered in a sequel paper.

More specifically, we will construct a one dimensional random field $(u(x))_{x \in \R}$ with remarkable multifractal and asymmetric (or skewed) properties: see equations \eqref{mainmultiproprerty} and \eqref{multiskewness} below. Though it is rather easy and classical to construct a field $u$ which satisfies the multifractal property \eqref{mainmultiproprerty} or the skewness property \eqref{multiskewness}, it is non trivial to construct a field $u$ which satisfies both properties. As will be explained in section \ref{Physicsmotivation}, both properties are essential in view of applications to turbulence.

On the mathematical side, the rigorous analysis of the field $u$ is also quite challenging and relies on a refined analysis of Gaussian multiplicative chaos measures (GMC, hereafter) and on the study of quite tedious integrals depending on real parameters. The theory of GMC measures has grown into an essential field of probability theory since the pioneering work of Kahane \cite{Kah85}: apart from turbulence, GMC measures are also widespread in the general field of conformal field theory (with applications to Liouville quantum gravity) and  in the field of finance (let us mention that modeling asset prices in a realistic way bears striking similarities with the topic of this paper, i.e. the modeling of the velocity field of a turbulent flow). See for instance the review \cite{RhoVar14} concerning this topic.

The organization of the paper is as follows. Next, we introduce the main notations of the paper and state the main results, i.e theorem \ref{th1secdef} and proposition \ref{propskewness} which state the important properties of the field $u$. In the next section, we present in detail the main motivations from turbulence that justify the construction of the field $u$. Then, in the following sections we proceed with the proofs.

\subsection{Notations and statement of the main results}
We consider a log-correlated stationary centered Gaussian field $\widehat{X}$ with covariance 
\begin{equation*}
\E[  \widehat{X}(x)\widehat{X}(0) ]= \ln_+\frac{L}{|x|}
\end{equation*}
for some fixed length scale $L>0$ where $\ln_+ x= \max( \ln x, 0 )$ for $x>0$. In the sequel, we will set $\frac{1}{|x|_+}=e^{ \ln_+\frac{1}{|x|}  }= \max( \frac{1}{|x|},1  )$.
We also consider a smooth regularization of $\widehat{X}$, call it $\widehat{X}_\epsilon$, with covariance structure satisfying
\begin{equation}\label{eq:DefCHat}
\forall x\ne 0, \; \; \; \; \; \widehat{C}_\epsilon(x)=\E[  \widehat{X}_\epsilon(x)\widehat{X}_\epsilon(0) ] \underset{\epsilon \to 0}{\rightarrow}  \ln_+\frac{L}{|x|},
\end{equation}
and set
\begin{equation}\label{eq:Defcepsilon}
c_\epsilon=\widehat{C}_\epsilon(0)=\E[\widehat{X}_\epsilon(x)^2].
\end{equation}
One can for instance choose $\widehat{X}_\epsilon= \frac{1}{\epsilon} (\widehat{X} \ast \theta (\frac{.}{\epsilon}))$  where $\theta$ is a smooth mollifier.

\medskip

We will focus in the rest of this paper on the following stochastic model of 1-d turbulent velocity field:

\begin{definition}[The stochastic field under study]\label{d.main}
 Recall that $c_\epsilon=\E[\widehat{X}_\epsilon(x)^2]$ \eqref{eq:Defcepsilon} and consider an independent white noise, call it $W$. We consider then a regularized field $u_\epsilon(x)$, $x\in \mathbb R$, defined by  
\begin{equation}\label{eq:SkMulProc}
u_\epsilon(x):=\int  \phi(x-y) X_\epsilon(y) e^{\gamma \widehat{X}_\epsilon(y) -\gamma^2 c_\epsilon}\,W(dy),
\end{equation}
where we have set
\begin{equation}\label{eq:DefX}
X_\epsilon(y):=\int k_\epsilon(y-z)e^{\gamma \widehat{X}_\epsilon(z) -\gamma^2 c_\epsilon}\,W(dz),
\end{equation}
 and the two following deterministic kernels:
 \begin{equation}\label{eq:DefphiIntro}
 \phi(x)=\varphi_L(x)\frac{1}{|x|^{\frac{1}{2}-H}},
 \end{equation}
with $\varphi_L(x)$ a characteristic cut-off function over the (large) fixed length scale $L$, that we assume without loss of generality to be even, and
 \begin{equation}\label{eq:DefkepsilonIntro}
k_\epsilon(x)=\frac{x}{|x|_\epsilon^{\frac{3}{2}}}1_{|x|\le L},
 \end{equation}
where  $|x|_\epsilon$ is a regularized norm of the vector $x$ over the (small) length scale $\epsilon$.
In the sequel, for the sake of clarity, we will only consider the case $L=1$ and will use $\varphi =\varphi_{1}$. This is no restriction as the general case can be dealt with similarly.

\end{definition}

It is obvious to check that the process $u_\epsilon$ is statistically homogeneous. The point is to determine whether the family of processes $(u_\epsilon)_\epsilon$ converges to a non trivial limit as $\epsilon\to 0$. For this, we will assume throughout this section that
\begin{equation}\label{gammasecdef}
2\gamma^2<1.
\end{equation}
The requirement \eqref{gammasecdef} is the optimal condition (\cite{RhoVar14}) to ensure   that almost surely  the random measures
 $$M_{2\gamma}^\epsilon(dy):=e^{2\gamma\hat{X}_\epsilon(y)-2\gamma^2 c_\epsilon}\,dy$$
converge weakly towards a random measure $M_{2\gamma}$ on $\R$. We also introduce the so-called Hurst index 
\begin{equation}\label{H}
H\in ]0,1[.
\end{equation}

Now we detail our main results. 
 \begin{theorem}\label{th1secdef}
 Assume \eqref{gammasecdef}+\eqref{H}. Then:
 \begin{enumerate}
 \item consider $q$ such that $0\leq q<\tfrac{1}{2\gamma^2}\wedge ( 1+\tfrac{H}{2\gamma^2})$.  We have for all $x\in\R$
 $$\sup_{\epsilon\in]0,1]}\E\Big[|u_\epsilon(x)|^q \Big] <+\infty.$$ 
 \item the marginals of the family $(u_\epsilon)_\epsilon$ converge in law as $\epsilon\to 0$ towards the marginals of some stationary centered  stochastic process $u$,  which is continuous and satisfies for  $0\leq q<\tfrac{1}{2\gamma^2}\wedge (1+\tfrac{H}{2\gamma^2})$
 \begin{equation}\label{mainmultiproprerty}
 \E[|u(x)-u(y)|^q]\sim C_q |x-y|^{\xi(q)}
 \end{equation}
 where 
 \begin{equation}\label{xi}
\xi(q)=(H+2\gamma^2)q-2\gamma^2q^2.
\end{equation}
 \end{enumerate}
 \end{theorem}
 This is the most general statement that we can claim under rather weak assumptions. The point that we want to improve is the continuity of sample paths of our process, and even Holder continuity. Thus we claim
 
 \begin{corollary}\label{coro2}
 Assume \eqref{gammasecdef}+\eqref{H}, and the further condition $H+(\sqrt{2}\gamma-1)^2>1$, then the process $u$ given by Theorem \ref{th1secdef} has 
  almost surely, continuous sample paths, which are even  locally  $\alpha$-H\"older for any $\alpha<H+ (\sqrt{2}\gamma-1)^2-1 $.
 \end{corollary}

\begin{rem}
In fact, we could certainly prove the almost sure uniform convergence over compact sets along subsequences of the family $(u_\epsilon)_\epsilon$ but we are more interested in the existence of the limiting process than the way it can be approximated. Furthermore we do not expect that the H\"older exponent we give above is  optimal.
\end{rem}
 
Coming back to our original motivations in turbulence, we want to make sure that the field, once the asymptotic limit $\epsilon\to 0$ has been taken, possesses moments of increments $u(x)-u(y)$ of order at least 3,  without the absolute value, and that they go to zero as a power-law of the distance $x-y$, reflecting the skewness of the field. The condition for existence of moments of order 3 is different from the one of the moments of increments, with absolute value, depicted in Theorem \ref{th1secdef}. Instead, we show the condition for existence of $ \E[(u(x)-u(y))^3]$ is 
\begin{equation}\label{mom3secdef}
\gamma^2<\frac{1}{8},
\end{equation}
independently on $H\in]0,1[$.
Our analysis of the moment of order 3, asymptotically in the limit $x-y\to 0^+$ (see Eq. \eqref{eq:EquivMom3SmallScales} for a complete expression), leads to
\begin{align*}
\E[(u(x)-u(y))^3]\propto (x-y)^{3H-12\gamma^2} \int_0^{\infty} f_H(h)\frac{1}{h^{\frac{1}{2}  +12 \gamma^2 }}dh,
\end{align*}
where is involved a special function $f_H(h)$ defined by a integral formula (see Eq. \eqref{deff_H}). Whereas we can compute the behaviors of $f_H(h)$ for small and large values of the argument $h$, and determine the conditions of the existence of this third order moment, overall the function $f_H$ is tedious to study. In particular, it turns out to be difficult to show that the integral entering in the former expression does not vanish. But a simple numerical estimation of the function $f_H$, presented in Annex \ref{ann:NumEstF_H}, suggests that the following holds
\begin{assumption} 
\begin{align}\label{conjf_H}
  \forall h>0,\mbox{ and }H\in ]0,1[/\{1/2\}, \;\;\; \left(\frac{1}{2}-H \right)f_{H}(h)>0. \end{align}
\end{assumption}
Thus, at a given parameter $H$, the function $f_H(h)$ does not change its sign, that would show that indeed the integral entering in the expression of the third order moment does not vanish.
 
Our next proposition relies on assumption \eqref{conjf_H} which seems challenging to prove rigorously (in spite of overwhelming numerical evidence that it is true):
 \begin{proposition}\label{propskewness}
 Assume \eqref{mom3secdef} and  \eqref{conjf_H}. Then 
 \begin{equation}\label{multiskewness}
 \E[(u(x)-u(y))^3]\build{\sim}_{x-y\to 0^+}^{} C'_3 (x-y)^{\xi(3)}
 \end{equation}
 for some non vanishing constant $C'_3$.
 \end{proposition}

 \begin{figure}[h]
\begin{center}
\includegraphics[width=13cm]{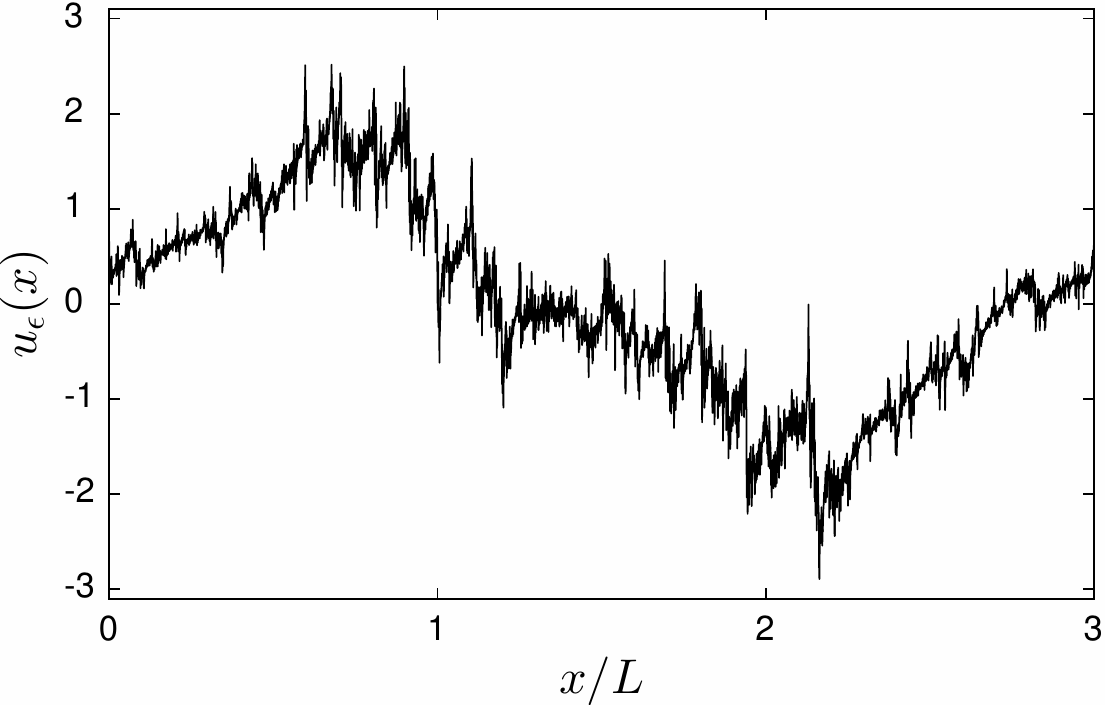}
\end{center}
\caption{\label{fig:ShowProc} An instance of the random process $u_\epsilon$ defined in \ref{eq:SkMulProc} over 3 cut-off length scale $L$. Here, we have renormalized the process by its standard deviation, and we have used $\gamma^2=0.025/4$ and $H=\frac{1}{3}+4\gamma^2$, at a given resolution scale $\epsilon$. See details in Annex \ref{ann:RandomSim}. }
\end{figure}

 \begin{rem}\label{rem:RemTurbuIntro}
As far as the modeling of turbulence is concerned, analysis of experimental measurements and numerical simulations of the Navier-Stokes equations give the universal value $4\gamma^2=0.025$ \cite{Fri95,CheCas12}. Given this value for $\gamma^2$, we are led to choose $H=\frac{1}{3}+4\gamma^2$, according to Eqs. \eqref{xi} and \eqref{multiskewness}, in order to fulfills the requirement of the  $4/5^{\mbox{\tiny th}}$-law of turbulence (see Section \ref{Physicsmotivation}) that states that $\xi(3)=1$. For this set of parameters, the multiplicative constant $C'_3$ entering in the proposition \ref{propskewness}, given the assumption \ref{conjf_H}, is strictly negative.

With these given values for the free parameters $\gamma$ and $H$, we represent in Fig. \ref{fig:ShowProc} an instance of the process $u_\epsilon$, at a given resolution scale $\epsilon$. See Annex \ref{ann:RandomSim} for details on the numerical simulation. We see by eyes that statistical laws are asymmetric at small scales. In particular, large negative values of increments are more probable than large positive ones.
\end{rem}

\subsection*{Acknowledgement}
The authors are partially supported by ANR grant \textsc{Liouville} ANR-15-CE40-0013. C.G. is furthermore partially supported by ERC grant LiKo 676999. L.C. thanks K. Gawedzki for fruitful discussions on singular integrals.

\section{Axiomatics of Kolmogorov's theory of turbulence, and design of the velocity field} 

\subsection{The phenomenological theory of Kolmogorov}\label{Physicsmotivation}

The statistical theory of incompressible, homogeneous and isotropic hydrodynamic turbulence is notoriously known to be a difficult matter. Making such a statement, as it can be already found in many classical textbooks \cite{Bat53,MonYag71,TenLum72,Fri95,Pop00}), would not surprise anyone since the underlying dynamics of viscous fluids is given by the Navier-Stokes equations, the study of which constitutes a difficult mathematical problem. The very link between the statistical approach and these dynamical equations is discussed in \cite{FoiMan01}, and is known in the physics literature as the Hopf's equation.

Based on natural and laboratory observations, the statistics of velocity fluctuations of fully developed turbulence are mostly understood in a phenomenological framework, given a limited set of free parameters, for which Kolmogorov made a series of key contributions \cite{Kol41,Kol62} (see the textbook of \cite{Fri95} for a extended presentation of this, and related numerous contributions of several authors).
 
The approach of Kolmogorov recasts the observed fluctuations of a fully developed turbulent velocity field, stirred at large scale by a stationary, say random, external divergence-free forcing vector field, into a consistent axiomatic framework \cite{Fri95}. To fix the ideas, call $u_i^\nu(x,t)_{i\in\{1,2,3\}, x\in\mathbb R^3, t\in\mathbb R}$ such a velocity field, and $\nu>0$ the kinematic viscosity of the fluid under interest. The time evolution of the velocity field is given by the incompressible Navier-Stokes equation (the density of the fluid is taken as unity) and reads
\begin{equation}\label{eq:NS}
\partial_t u^\nu_i + u^\nu_j\partial_ju^\nu_i=-\partial_ip^\nu+\nu\Delta u^\nu_i + f_i
\end{equation}
where $p$ is the pressure field, determined by the additional  incompressibility condition $\partial_iu^\nu_i=0$, and $f$ a divergence forcing field, smooth and typically correlated over a large spatial scale $L$, called the integral length scale in turbulence literature. This large scale $L$ is schematically the scale at which energy is injected into the flow, and is independent of the viscosity $\nu$.

The phenomenology of Kolmogorov can be decomposed in terms of 3 axioms, that remain for the most part, as far as we know, unrelated to the Navier-Stokes equations. For simplicity, we will present them in a uni-dimensional context. Take for instance the velocity component along a given axis in the laboratory reference frame, say $x$, and consider its spatial distribution along that very same axis. This field can be measured experimentally in wind tunnels or in jets, once the so-called streamwise velocity component is interpreted in a spatial context using the Taylor's frozen hypothesis \cite{Fri95}. Henceforth, we work in the 1-dimensional space, and we call the respective velocity field $u^\nu(x,t)_{ x\in\mathbb R, t\in\mathbb R}$.  The axioms read
\begin{itemize}
\item \textit{(concerning the velocity variance)} When forced by the divergence-free vector field $f$ entering in the Navier-Stokes equations (Eq. \ref{eq:NS}), the velocity field reaches a statistically stationary regime, in which the variance of any velocity components remains finite and becomes independent on viscosity $\nu$ when $\nu$ gets smaller. In particular, as far as the velocity component we are interested in is concerned, we write
\begin{equation}\label{eq:Ax1}
 \sigma^2 = \lim_{\nu\to 0}\lim_{t\to \infty}\E \left[\left(u^\nu(x,t)\right)^2\right] < \infty.
\end{equation}

\item \textit{(concerning the asymptotics of the mean dissipation)}  To ensure the finiteness and $\nu$-independence of the variance, the flow will self-organize to dissipate all the injected energy at a $\nu$-independent rate. This axiom reads more precisely
\begin{equation}\label{eq:Ax2}
0<\lim_{\nu\to 0}\lim_{t\to \infty} \nu\E \left[|\partial_x u^\nu(x,t)|^2\right]< \infty.
\end{equation}
In other words, as $\nu$ gets smaller, this average \textit{viscous} dissipation becomes independent on $\nu$ itself. It is known in the literature as the \textit{dissipative anomaly} \cite{EyiSre06}. Its actual value can then only be given by the statistical properties of the flow at large scale $L$. From a dimensional point of view, there is no other choice than $\sigma^3/L$.

\item \textit{(concerning the asymptotical non-differential nature of the velocity field)} As depicted by the second axiom (Eq. \ref{eq:Ax2}), as viscosity vanishes, the variance of the spatial gradients diverges (as $1/\nu$). In other words, as $\nu\to 0$, the velocity field remains bounded but is rough (i.e. non differentiable). In particular,  for $q\in \mathbb N$, the respective structure functions behave as
\begin{equation}\label{eq:Ax3}
\lim_{\nu\to 0}\lim_{t\to \infty} \E \left[\left(u^\nu(x,t)-u^\nu(y,t)\right)^q\right]\build{\sim}_{x-y\to 0^+}^{} D_q \sigma^q\left(\frac{x-y}{L}\right)^{\xi_q},
\end{equation}
where $D_q$ and $\xi_q$ are universal functions of the order $q$, universal in the sense that they are independent on both characteristic scale and amplitude of the forcing term, and on viscosity $\nu$. Here, the constants $D_q$, up to nondimensionalization  using $\sigma$ and $L$,  are related to the constants $C_q$ (when $q$ is even) and $C'_3$ (when $q=3$) entering respectively in theorem \ref{th1secdef} and proposition \ref{propskewness}. When looking at experimental data, we can estimate that, in good approximation, $\xi_2\approx 2/3$. This is called the $2/3^{\mbox{\tiny th}}$-law of turbulence, which is not based on a rigorous derivation assuming axioms 1 (Eq. \ref{eq:Ax1}) and 2 (Eq. \ref{eq:Ax2}). It corresponds in a equivalent formulation in Fourier space to the power-law decay of the velocity power-spectrum with an exponent $5/3$, and says that the H$\ddot{\mbox{o}}$lder exponent of velocity is close to $1/3$ in a statistically averaged sense. Furthermore, assuming the first two axioms (Eq. \ref{eq:Ax1} and \ref{eq:Ax2}), using the stationary solution of the so-called K\'arm\'an and Howarth equation \cite{Kol41,Fri95}, it can be shown rigorously that $\xi_3=1$, and that $D_3$ is a universal constant, strictly negative, and related to the ratio of the average viscous dissipation by its naive dimensional estimation $\sigma^3/L$ (multiplied by the factor $-4/5$). This is called the $4/5^{\mbox{\tiny th}}$-law of turbulence \cite{Fri95}. Processing experimental data beyond the second and third order moments suggests strongly that $\xi_q$ is a non linear and concave function of the order $q$. This is known as the intermittency, i.e. multifractal, phenomenon. We note, as a definition, $\gamma^2 = -\frac{1}{4}(\partial^2\xi_q/\partial q^2)_{q=0}$ (consistently with \eqref{mainmultiproprerty}) the intermittency coefficient. It is observed universal, i.e. for any Reynolds numbers and forcing conditions, and it has been measured that $4\gamma^2=0.025$ (see for instance \cite{Fri95,CheCas12} and references therein).
\end{itemize}
 
\subsection{The underlying fractional Gaussian field}

To go further in this statistical picture of turbulence, we could wonder whether it is possible to give a rigorous meaning of this ensemble of three axioms. We are thus asking whether we can build up a 1d-velocity field that mimics the fluctuations of fully developed turbulence, using, as it is classically done in a Wiener chaos expansion, Wiener integrals \cite{Nua00}. A proposition of such a stochastic representation of turbulence was firstly made by Kolmogorov, and formalized by Mandelbrot and van Ness \cite{ManVan68} in the more general class of Gaussian fractional Brownian motions. Call $u_\epsilon^{\mbox{\tiny g}}(x)_{x\in \mathbb R}$ such a Gaussian process, defined by
\begin{equation}\label{eq:ugepsilon}
u_\epsilon^{\mbox{\tiny g}}(x)=\int \phi_\epsilon(x-y)\,W(dy),
\end{equation} 
where $\phi_\epsilon$ is a deterministic kernel given by
 \begin{equation}\label{eq:KerFbF}
\phi_\epsilon(x)=\varphi_L(x)\frac{1}{|x|_\epsilon^{\frac{1}{2}-H}},
 \end{equation}
and $|x|_\epsilon$ a regularized norm of the vector $x$ over the (small) length scale $\epsilon$, and $\varphi_L(x)$ a characteristic cut-off function over the (large) length scale $L$. Adapting the arguments of Ref. \cite{ManVan68}, as it is done in Refs. \cite{RobVar08} and recalled in \cite{PerGar16}, it can be shown that the Gaussian process $u_\epsilon^{\mbox{\tiny g}}(x)$ (Eq. \ref{eq:ugepsilon}) converges when $\epsilon\to 0$, whatever the regularizing mechanism (entering in the very definition of the regularized norm $|x|_\epsilon$) and for $H\in]0,1[/\{1/2\}$, towards a finite variance Gaussian process $u^{\mbox{\tiny g}}(x)$, which is non-differentiable. Its structure functions are given by
\begin{equation}\label{eq:StructUg}
 \E \left[\left(u^{\mbox{\tiny g}}(x)-u^{\mbox{\tiny g}}(y)\right)^q\right]\build{\sim}_{x-y\to 0^+}^{} D^{\mbox{\tiny g}}_q \sigma^q\left(\frac{x-y}{L}\right)^{qH},
\end{equation}
where $D^{\mbox{\tiny g}}_q$ are universal constants, in the sense we have defined before. Moreover, $D^{\mbox{\tiny g}}_{2q+1}=0$. Going back to the physics of turbulence, considering the particular case $H=1/3$, we see that this Gaussian field fulfills axiom 1 (Eq. \ref{eq:Ax1}), the regularizing scale $\epsilon$ can be eventually chosen with the appropriate dependence on $\nu$ in order to fulfill axiom 2 (Eq. \ref{eq:Ax2}), but fails at reproducing both the $4/5^{\mbox{\tiny th}}$-law of turbulence (i.e. $C^{\mbox{\tiny g}}_3=0$) and the nonlinear behavior of the spectrum of exponents $\xi_q$ related to its non-vanishing curvature at the origin ($\gamma\ne 0$). Nonetheless, a Gaussian process as a underlying random field is an appealing starting point since it allows to reproduce the first two axioms and the $2/3^{\mbox{\tiny th}}$-law of turbulence. The purpose of this article is to show how to properly disturb this Gaussian field in order to reproduce the missing key behaviors of turbulence, which are the $4/5^{\mbox{\tiny th}}$-law and the intermittency phenomenon.

\subsection{Introduction of the Gaussian multiplicative chaos}

To do so, let us first focus on the intermittency phenomenon. One of the simplest way to understand  and model this intrinsically  non Gaussian nature of turbulence and its related non-vanishing intermittency parameter $\gamma\ne 0$, is to consider the exponential of a logarithmically correlated Gaussian field $X(x)$, as it was proposed by Mandelbrot \cite{Man72}. In this spirit, the theory of multiplicative chaos \cite{Kah85,RhoVar14} gives a rigorous meaning to the random measure ``$e^{\gamma X}$" applying a proper regularizing procedure and taking the limit.

Remark that the Gaussian process $X$ is assumed to be logarithmically correlated, in particular the variance has to diverge. One could then wonder what meaning can be given to the exponential of it. This is properly understood in the framework of multiplicative chaos \cite{Kah85,RhoVar14}. Disturbing the underlying Gaussian field $u^{\mbox{\tiny g}}(x)$ previously described would be naturally done while multiplying the Gaussian white measure $W$ entering in its definition (c.f. \ref{eq:ugepsilon}) by this singular measure ``$e^{\gamma X}$". Thus, we are asking for a proper meaning of not only the exponential of the logarithmically correlated Gaussian process $X$, but also the meaning of this multiplicative chaos multiplied by a distributional white noise $W$. Preliminary mathematical study of this was done in Ref. \cite{RobVar08}, also in the context of fully developed turbulence. It turns out that giving a meaning to the random distribution  ``$e^{\gamma X(x)}W(dx)$" is not obvious, besides the trivial case of taking $X$ and $W$ independent. It is easy to see that if indeed the multiplicative chaos and the white measure are taken independent, this will lead to a vanishing third order moment of velocity increment, and thus to the impossibility of reproducing the $4/5^{\mbox{\tiny th}}$-law of turbulence.  The difficulty in defining  ``$e^{\gamma X(x)}W(dx)$" as a well defined random measure lies in the necessity, in particular for turbulence modeling purposes, to consider the two fields $X$ and $W$ being correlated. Let us also keep in mind that at the end, we would like furthermore to apply a linear operation on this measure with a kernel $\phi$ (Eq. \ref{eq:KerFbF}) that becomes singular in the range of interesting values $H<1/2$.

In this context, considering, instead of the product of the two random distributions $e^{\gamma X(x)}$ and $W(dx)$, the product of the distribution $e^{\gamma X(x)}$ by a finite-covariance field, say $\omega(x)dx$, appears to be a natural way to properly define such a measure. This has been studied in Ref. \cite{BacDuv12} for financial applications. It turns out that their proposition, that is to take for the Gaussian field $\omega(x)$ a fractional Gaussian noise as considered in \cite{ManVan68} (see Ref. \cite{BacDuv12} for details) and a given cross-covariance structure of the fields $\omega(x)$ and $X(x)$, allows to reproduce a non-vanishing third-order structure function when $H>1/2$, making such a process \textit{skewed}. Unfortunately for applications in turbulence theory, the interesting case $H<1/2$ is not included in their results. Further theoretical works in this direction \cite{AbrCha09,Per13}, taking as a simplified framework the case of independent fields $\omega(x)$ and $X(x)$, indeed show that the predicted spectrum of exponents $\xi_q$ becomes independent on $H$ when $H<1/2$. This is related to the pathological behavior of the fractional Gaussian noise at a high level of roughness $H<1/2$.  

We see that constructing a uni-dimensional random field that fulfills the three axioms of Kolmogorov, including intermittency $\gamma\ne 0$, the skewness phenomenon (i.e. $D_3\ne 0$) and $\xi_2$ close to 2/3, is far from being easy. We thus need to rethink the building block of the Wiener integrals that we want to make use of. Elaborating on the propositions made in Ref. \cite{RobVar08}, it was proposed in Ref. \cite{CheRob10}, to include in this picture some aspects of the Euler equations, and more precisely some aspects of the vorticity stretching mechanism. The main output of this work is the proposition of a homogeneous, isotropic, incompressible (i.e. divergence-free) random vector field, based on a underlying fractional Gaussian field structure and a matrix multiplicative chaos (developed in Ref. \cite{CheRho13}) that is shown, numerically, to be realistic of a fully developed turbulent flow. In other words, as far as we could go in a numerical simulation, the proposed velocity field of Ref. \cite{CheRho13} is consistent with the axiomatic approach of Kolmogorov. Unfortunately, at this stage, an exact derivation of its statistical properties, and their asymptotical behavior when the resolution scale $\epsilon$ goes to zero, seems to be out of range. The difficulty in obtaining exact properties is linked to the existence of correlations between the matrix multiplicative chaos and the underlying vector white noise entering in the construction. A further theoretical analysis of simplified versions of this random vector field, assuming for instance independence of the matrix multiplicative chaos and the underlying Gaussian white vector, or performing a perturbative expansion in the small parameter $\gamma$, was proposed in Ref. \cite{PerGar16}. This study confirmed the realism of the random vector field using massive numerical simulations, and illustrated some of the mechanisms at the root of a non-vanishing third order structure function. 

\subsection{An intermediate non canonical uni-dimensional field}

Given the observed (numerically) realism of the vector field proposed in Refs. \cite{CheRob10,PerGar16}, and in front of the difficulty of defining rigorously the product of a (matrix) multiplicative chaos by a (vector) white noise with a given inter-dependence, we make hereafter the choice to work with a uni-dimensional ersatz. This allows us to push forward our understanding of such singular measures, and to obtain exact expressions of their statistical properties. A natural choice of this ersatz, in the spirit of the stochastic structure of the formerly described vector field, is to consider the following intermediate field
\begin{equation}\label{eq:uintepsilon}
u_\epsilon^{\mbox{\tiny int}}(x)=\int\phi(x-y)e^{\gamma \tilde{X}_\epsilon(y)-\gamma^2\E(\tilde{X}_\epsilon^2)}\,W(dy),
\end{equation} 
where enters, compared to the Gaussian version of the field \eqref{eq:ugepsilon} ($\phi$ being the pointwise limit of the deterministic kernel $\phi_\epsilon$ defined in \eqref{eq:KerFbF} when $\epsilon\to 0$), an additional multiplicative chaos obtained while exponentiating a Gaussian field $\tilde{X}_\epsilon$. It is defined as the following stochastic integral
\begin{equation}\label{eq:DefXtilde}
\tilde{X}_\epsilon(y):=\int k_\epsilon(y-z)\,W(dz),
\end{equation}
where enters a deterministic kernel given by
 \begin{equation}\label{eq:KerLCF}
k_\epsilon(x)=\frac{x}{|x|_\epsilon^{\frac{3}{2}}}1_{|x|\le L},
 \end{equation}
that ensures to $\tilde{X}_\epsilon$ a logarithmic correlation structure in the limit $\epsilon\to 0$. It is a standard construction in the framework of the GMC \cite{RhoVar14}, that is made to give intermittent, i.e. multifractal, corrections to the underlying Gaussian structure \eqref{eq:ugepsilon}. Similar types of fields, such as \eqref{eq:uintepsilon}, were considered in Ref. \cite{RobVar08}, but here, and it is original at this stage, the kernel entering the construction \eqref{eq:KerLCF} is odd, as suggested by the tensorial kernels entering in the definition of the vector fields of Refs. \cite{CheRob10,PerGar16}. Notice here, and it is crucial to allow for a possible non vanishing third order moment of velocity increments, that the very same instance of the white noise $W$ enters in both the first layer of $u_\epsilon^{\mbox{\tiny int}}$ \eqref{eq:uintepsilon}, and the field $\tilde{X}_\epsilon$ \eqref{eq:DefXtilde}, imposing a complex internal correlated structure. In particular, we have $\E[\tilde{X}_\epsilon(x)W(dy)] = k_\epsilon(x-y)dy$. As it is shown in Refs. \cite{CheGar15,CheHDR15}, it becomes necessary here to specify the regularization procedure. For instance, take $|.|_\epsilon =\frac{1}{\epsilon} (|.| \ast \theta (\frac{.}{\epsilon}))$  where $\theta$ is a smooth mollifier, and remark that $|\epsilon .|_\epsilon =(|.| \ast \theta (.)) = |.|_1$. Consider then the rescaled (by $\epsilon$) quantities
$$ r_\epsilon(h) = \sqrt{\epsilon}k_\epsilon(\epsilon h) =\frac{h}{|h|_1^{\frac{3}{2}}}1_{|h|\le L/\epsilon} \build{\rightarrow}_{\epsilon\to 0}^{} r(h)=\frac{h}{|h|_1^{\frac{3}{2}}} \,\,\mbox{ for } h\in \mathbb R,$$
and
$$R_\epsilon(h)= \E \left[\tilde{X}_\epsilon(0)\tilde{X}_\epsilon(\epsilon h)\right] - \E(\tilde{X}_\epsilon^2)  \build{\rightarrow}_{\epsilon\to 0}^{} R(h)= \int r(x) \left[ r(x+h)-r(x)\right]\,\,\mbox{ for } h\in \mathbb R. $$
Using the Gaussian integration by parts \cite{Nua00}, it is easy to show that the field $u_\epsilon^{\mbox{\tiny int}}$ \eqref{eq:uintepsilon} is centered, and we can then obtain \cite{CheGar15,CheHDR15}
$$ \E \left[\left(u_\epsilon^{\mbox{\tiny int}}\right)^2\right] \build{\rightarrow}_{\epsilon\to 0}^{}  \E \left[\left(u^{\mbox{\tiny int}}\right)^2\right] = \E \left[\left(u^{\mbox{\tiny g}}\right)^2\right] \left[ 1-\gamma^2 \int r^2(h)e^{\gamma^2 R(h)}dh\right].$$
This former expression shows that asymptotically, in the limit $\epsilon\to 0$, the variance of the field $u^{\mbox{\tiny int}}$ depends on how its approximation at the scale $\epsilon$ has been made, and in particular on the precise choice of the mollifier $\theta$. From the physical point of view, although the variance of this field is finite and not vanishing, it contradicts somehow the first axiom of the phenomenology of Kolmogorov \eqref{eq:Ax1} since it is not expected that the mechanisms at play in the viscous dissipation would contribute: only the forcing field and boundary conditions should determine the variance. In this sense, we will say that the convergence of the field (as $\epsilon\to 0$) towards its asymptotical form is not canonical: it keeps track of the choice that has been made to regularize the field at a given scale $\epsilon>0$.

Furthermore, even if the field has a non vanishing third order moment of increments at a finite $\epsilon>0$, it can be demonstrated  \cite{CheGar15,CheHDR15} that it looses this property when $\epsilon\to 0$, i.e. for any $x\ne y$, we have asymptotically
\begin{equation}\label{eq:LossM3}
\E \left[\left(u_\epsilon^{\mbox{\tiny int}}(x)-u_\epsilon^{\mbox{\tiny int}}(y)\right)^3\right]\build{\rightarrow}_{\epsilon\to 0}^{} 0, 
\end{equation}
showing in a definitive way that the field  $u^{\mbox{\tiny int}}$ is unable to reproduce in a realistic manner the set of axioms presented in Section \ref{Physicsmotivation}.

\subsection{Present approach}

Although the intermediate field $u_\epsilon^{\mbox{\tiny int}}$ \eqref{eq:uintepsilon} that we considered in the former section exhibits a non canonical way of convergence when $\epsilon \to 0$ (with a loss of the third order moment of increments \eqref{eq:LossM3}), it shows that giving a meaning to ``$e^{\gamma X(x)}W(dx)$" is non trivial. In particular, if is assumed a odd-correlation structure between the logarithmically correlated field $X$ and the white noise $W$ as it is done in \eqref{eq:uintepsilon}, the statistical laws of the asymptotical field $u^{\mbox{\tiny int}}$ keeps track of the regularization procedure that what used to define it. Let us then keep in mind that a simple and canonical way to define the random measure ``$e^{\gamma X(x)}W(dx)$" is to consider the fields $X$ and $W$ as being independent. We now need to propose another way to introduce a correlated internal structure to the field in order to model asymmetric probability laws for the velocity increments.

To do so, let us elaborate on the former field $u_\epsilon^{\mbox{\tiny int}}$ \eqref{eq:uintepsilon} and remark, at a given $\epsilon>0$, that it can be developed in powers the intermittency parameter $\gamma$ such as to obtain
$$ u_\epsilon^{\mbox{\tiny int}}(x)=u_\epsilon^{\mbox{\tiny g}}(x)+\gamma\int\phi(x-y)\tilde{X}_\epsilon(y)\,W(dy)+o_\epsilon(\gamma),$$
where $o_\epsilon(\gamma)$ stand for a random field made up of the higher order terms of the development of $ u_\epsilon^{\mbox{\tiny int}}$ in powers of $\gamma$. We see that the second field entering in the former development (proportional to $\gamma$) coincides exactly with our present field 
$ u_\epsilon$ \eqref{eq:SkMulProc} when its parameter $\gamma$ is set to 0. Call it  $ (u_\epsilon)_{\gamma=0}$. It turns out that this field is a well-defined random object in the limit $\epsilon\to 0$, with in particular a finite variance (see devoted theoretical material leading to \eqref{eq:AsymptVarMonofractal}) that does not depend on the regularization procedure, and furthermore exhibits a non vanishing third order moment of velocity increments (take $\gamma=0$ in \eqref{eq:Mom3IncrMulti} and see \eqref{eq:EquivMom3SmallScalesGamma0} for the behavior at small scales). Thus, the naive development in powers of $\gamma$ that we performed exhibits a field that appears a to be good candidate to obtain a non vanishing third order moment of increments.

It remains to introduce in $ (u_\epsilon)_{\gamma=0}$, as a final layer, the intermittency corrections. This is done while considering an independent multiplicative chaos, and thus an independent logarithmically correlated process $\widehat{X}_\epsilon$,  and replacing the white noise $W$ entering in its expression by $e^{\widehat{X}_\epsilon} W$. As it is shown in \cite{DomRho11}, the multiplicative chaos has to be also introduced in a similar fashion in the field $\tilde{X}_\epsilon$ that enters the construction, and thus replacing $\tilde{X}_\epsilon$ by the field $X_\epsilon$ defined in \eqref{eq:DefX}. This final step is necessary to guaranty the power-law behaviors announced in Theorem \ref{th1secdef} and Proposition \ref{propskewness}. Doing so, we end up with the field $u_\epsilon$ \eqref{eq:SkMulProc} that we propose to study.

\section{Proof of the main results}

Now, we address the main results of the introduction by studying the statistical properties of $u_\epsilon$  \eqref{eq:SkMulProc}.

\subsection{Study of the average}

Because the field is statistically homogeneous, consider only the average of the field at the position $x=0$. Recall that $\phi$ is a even function of its argument, since we assumed $\varphi$ to be even itself. We then have
\begin{align*}
  \E[  u_\epsilon(0)  |   \widehat{X}_\epsilon  ] &  =     
\E\Big[\int \phi(-y)   \left ( \int k_\epsilon(y-z)e^{\gamma \widehat{X}_\epsilon(z) -\gamma^2 c_\epsilon}\,W(dz) \right )  e^{\gamma \widehat{X}_\epsilon(y) -\gamma^2 c_\epsilon}\,     W(dy)  \Big|   \widehat{X}_\epsilon  \Big]  \\
&  =   \iint  \phi(y)   k_\epsilon(y-z)e^{\gamma \widehat{X}_\epsilon(z) -\gamma^2 c_\epsilon}  e^{\gamma \widehat{X}_\epsilon(y) -\gamma^2 c_\epsilon}\,   \E[  W(dy)W(dz)]    \\
&  =  k_\epsilon(0) \int   \phi(y)  e^{2 \gamma \widehat{X}_\epsilon(y) -2\gamma^2 c_\epsilon} \, dy \\   &=0,
\end{align*}
since $k_\epsilon(0)=0$.  Thus, the random field is centered,  i.e. $\E[  u_\epsilon(0) ]=0$.

\subsection{Study of the variance}


Consider first the case $\gamma\ne 0$. Notice that   
\begin{equation*}
\E \left[X_\epsilon(y) W(dz)   \Big|   \widehat{X}_\epsilon  \right] = k_\epsilon(y-z)  e^{ \gamma \widehat{X}_\epsilon(y) -\gamma^2 c_\epsilon} dz,
\end{equation*}
and recall that $k_\epsilon$ is a odd function of its argument, in particular $k_\epsilon(0)=0$. Therefore, we have by integration by parts, recall that $\phi$ is even,
\begin{align*}
&\E \left[ u_\epsilon(0)^2  \Big|   \widehat{X}_\epsilon  \right] \\
&=\iint \phi(-y)\phi(-z)\E \left[X_\epsilon(y)X_\epsilon(z)\,W(dy)W(dz)  \Big|   \widehat{X}_\epsilon  \right] e^{ \gamma \widehat{X}_\epsilon(y) -\gamma^2 c_\epsilon}  e^{ \gamma \widehat{X}_\epsilon(z) -\gamma^2 c_\epsilon}  dydz \\
&=\int \phi^2(y)\E \left[X_\epsilon(y)^2  \Big|   \widehat{X}_\epsilon  \right] e^{2 \gamma \widehat{X}_\epsilon(y) - 2 \gamma^2 c_\epsilon} dy-\iint \phi(y)\phi(z)k_\epsilon^2(y-z)\, e^{ 2\gamma \left[\widehat{X}_\epsilon(y)+\widehat{X}_\epsilon(z) \right] -4\gamma^2 c_\epsilon}   dydz,
\end{align*}
where a computation on the white noise leads to
\begin{equation*}
\E \left[X_\epsilon(y)^2  \Big|   \widehat{X}_\epsilon  \right] = \int k^2_\epsilon(y-z)  e^{2 \gamma \widehat{X}_\epsilon(z)- 2 \gamma^2 c_\epsilon } dz.
\end{equation*}
Using this former expression of the conditional variance, we arrive at
\begin{align}\label{eq:ConvEpsVarCond}
\E \left[ u_\epsilon(0)^2  \Big|   \widehat{X}_\epsilon  \right] =\iint \left[\phi^2(y)-\phi(y)\phi(z)\right]k_\epsilon^2(y-z)\, e^{ 2\gamma \left[\widehat{X}_\epsilon(y)+\widehat{X}_\epsilon(z) \right] -4\gamma^2 c_\epsilon}   dydz.
\end{align}
Averaging on the field $\widehat{X}_\epsilon$, and symmetrizing the expression, we get
\begin{align*}
\E \left[ u_\epsilon(0)^2 \right] =\frac{1}{2}\iint \left[\phi(y)-\phi(z)\right]^2 k_\epsilon^2(y-z)\, e^{4\gamma^2\E[\widehat{X}_\epsilon(y)\widehat{X}_\epsilon(z)]}  dydz,
\end{align*}
that eventually converges towards a finite limit 
\begin{align}\label{eq:ConvVarSym}
\E \left[ u^2 \right] = \lim_{\epsilon\to 0}\E \left[ u_\epsilon(0)^2 \right] =\frac{1}{2}\iint \left[\phi(y)-\phi(z)\right]^2\frac{1}{|y-z|^{1+4 \gamma^2}} 1_{|y-z| \leq 1} dydz.
\end{align}
We indeed show in Section \ref{ss.mom} that the double integral entering in \eqref{eq:ConvVarSym} exists and is finite, for a certain range of values of $H$ and $\gamma$ (i.e. $\gamma^2 < \frac 1 2 H$). Thus by dominated convergence, it shows that the variance of the process $u_\epsilon$ \eqref{eq:SkMulProc} converges towards a finite, non vanishing and positive value which is independent of the regularization mechanism that we have chosen.

In order to see more clearly the underlying phenomena (and their cancellations) that take place behind this convergence, we propose to present a more straightforward way to show the convergence of the variance. This will also allow us to define key quantities that will be entering in the computation of the skewness (Section \ref{Sec:Skewness}).

Going back to \eqref{eq:ConvEpsVarCond}, by taking the expectation with respect to $\widehat{X}_\epsilon$, making the change of variable $h=y-z$ and integrating over $h$ and $y$, this yields 
\begin{equation*}
\E \left[ u_\epsilon(0)^2 \right]
=(\phi\star\phi)(0)  A_\epsilon -2\int_0^{\infty} (\phi\star\phi)(h)k_\epsilon^2(h)e^{4\gamma^2\E[\widehat{X}_\epsilon(h)\widehat{X}_\epsilon(0)]}\, dh,
\end{equation*}
where
\begin{equation*}
A_\epsilon=\int k_\epsilon^2(u)e^{4\gamma^2\E[\widehat{X}_\epsilon(u)\widehat{X}_\epsilon(0)]}\,du,
\end{equation*}
and where we use the following notation 
\begin{equation*}
( \phi\star\phi) (h) = \int \phi(x) \phi(x+h)  dx.
\end{equation*}
Note that $\phi\star\phi$ differs from the standard convolution.
\medskip

Define $K_\epsilon(h) = -\int_{h}^{\infty}k_\epsilon^2(x)e^{4\gamma^2\E[\widehat{X}_\epsilon(x)\widehat{X}_\epsilon(0)]}\, dx$ such that $K_\epsilon '(h) = k_\epsilon^2(h)e^{4\gamma^2\E[\widehat{X}_\epsilon(h)\widehat{X}_\epsilon(0)]}$ and $2K_\epsilon(0) = - A_\epsilon$. Remark that pointwise
\begin{equation*}
\lim_{\epsilon \to 0}K_\epsilon(h) =\frac{1}{4\gamma^2}\big(1-|h|^{-4\gamma^2}\big) 1_{|h| \leq 1}.
\end{equation*}
An integration by parts gives
\begin{align*}
\E u_\epsilon(0)^2&=(\phi\star\phi)(0) A_\epsilon -2\left[-(\phi\star\phi)(0)K_\epsilon(0)-\int_0^{\infty} (\phi\star\phi)'(h)K_\epsilon(h)\, dh\right]\\
&=2\int_0^{\infty} (\phi\star\phi)'(h)K_\epsilon(h)\, dh.
\end{align*}
Thus the variance converges by dominated convergence towards a finite value with 
\begin{align}\label{eq:AsymptVarMultifractal}
\E u(0)^2=\lim_{\epsilon \to 0}\E u_\epsilon(0)^2 = \frac{1}{2\gamma^2}\int_0^{1} (\phi\star\phi)'(h)\big(1-|h|^{-4\gamma^2}\big)\, dh.
\end{align}

In order to make sense of the asymptotical form of the variance \eqref{eq:AsymptVarMultifractal}, one has to check the integrability in the neighborhood of the origin. This is the subject of Lemma \ref{lem:Equiv2}. We indeed show that the function $(\phi\star\phi)(h)$ is continuously differentiable over $\mathbb R$ when $H\in ]1/2,1[$. For  $H\in ]0,1/2[$, the function $(\phi\star\phi)(h)$ is not differentiable at the origin. Consistently, for $H\in ]0,1[$,  the derivative $(\phi\star\phi)'(h)$ behaves at the origin as $|h|^{2H-1}$, a behavior which is integrable at the origin. This gives a meaning to the asymptotical variance \eqref{eq:AsymptVarMultifractal} as long as $1-2H+4 \gamma^2<1$. (i.e. $\gamma^2 < \frac 1 2 H$). 

 \begin{lemma}\label{lem:Equiv2}
For $H\in ]0,1[/\{1/2\}$  and $|h|>0$, the function $(\phi\star\phi)(h)$ is differentiable and its derivative is given by
\begin{align}\label{eq:ExpDerivEquiv2}
(\phi\star\phi)'(h) = &\int\varphi(x)\varphi'(x+h)\frac{1}{|x|^{\frac{1}{2}-H}}\frac{1}{|x+h|^{\frac{1}{2}-H}} dx\notag \\
&+ (H-1/2)\,\mathrm{ P.V. }\int\varphi(x)\varphi(x+h)\frac{1}{|x|^{\frac{1}{2}-H}}\frac{x+h}{|x+h|^{\frac{5}{2}-H}} dx,
\end{align}
where we have defined the principal value integral (P.V.) that can be written using a convergent integral:
$$\mathrm{ P.V. }\int\varphi(x)\varphi(x+h)\frac{1}{|x|^{\frac{1}{2}-H}}\frac{x+h}{|x+h|^{\frac{5}{2}-H}}dx  =\int_0^{\infty}\frac{\varphi(x)}{x^{\frac{3}{2}-H}}\left[\frac{\varphi(x-h)}{|x-h|^{\frac{1}{2}-H}}-\frac{\varphi(x+h)}{|x+h|^{\frac{1}{2}-H}}\right]dx. $$
Furthermore, for $H>1/2$, $(\phi\star\phi)'(h)$ is continuous, bounded over $\R$ and $(\phi\star\phi)'(0)=0$,  and we have for $H\in ]0,1[/\{1/2\}$  the following equivalent at the origin
 $$(\phi\star\phi)'(h)\build{\sim}_{h\to 0^+}^{}(H-1/2)\varphi^2(0)\sign(h)|h|^{2H-1}\,\mathrm{ P.V. }\int\frac{1}{|x|^{\frac{1}{2}-H}}\frac{x+1}{|x+1|^{\frac{5}{2}-H}} dx.$$

 \end{lemma}
 
 \noindent {\it Proof.} To prove the expression for the derivative \eqref{eq:ExpDerivEquiv2}, regularize the singularity and pass to the limit. To get the equivalent, rescale the dummy integration variable by $|h|$ in the second term of the RHS of \eqref{eq:ExpDerivEquiv2} and take the limit. Remark that this equivalent is also correct for $H>1/2$, since the first term of the RHS of \eqref{eq:ExpDerivEquiv2} behaves as $h$ at the origin (using the fact that $\varphi$ is even), i.e.
\begin{align*}
\int\varphi(x)&\varphi'(x+h)\frac{1}{|x|^{\frac{1}{2}-H}}\frac{1}{|x+h|^{\frac{1}{2}-H}} dx\\
&\build{\sim}_{h\to 0}^{}h \left[\int\varphi(x)\varphi''(x)\frac{1}{|x|^{1-2H}} dx+ (H-1/2)\,\mathrm{ P.V. }\int\varphi(x)\varphi'(x)\frac{x}{|x|^{3-2H}} dx\right],
\end{align*}
and so tends to 0 when $h\to 0$ faster than the second term.

\qed

\begin{rem}
In the case $\gamma=0$, one gets the following formula by similar computations
\begin{equation}\label{eq:AsymptVarMonofractal}
\E u(0)^2=\lim_{\epsilon \to 0}\E u_\epsilon(0)^2 = 2 \int_0^{1} (\phi\star\phi)'(h)\ln h \, dh.
\end{equation}
Notice that \eqref{eq:AsymptVarMonofractal} can also be obtained by taking the limit $\gamma \to 0$ in \eqref{eq:AsymptVarMultifractal}.
\end{rem}

\subsection{Study of the variance of increments}\label{Sec:VarVeloIncr}

Once again, consider first the case $\gamma\ne 0$, we will treat the case $\gamma=0$ as a remark at the end of this Section. Define the increments as
$$\delta_\ell u_\epsilon(x) = u_\epsilon (x+\ell/2)-u_\epsilon (x-\ell/2)=\int _{\mathbb R}\Phi_\ell(x-y)  X_\epsilon(y) e^{\gamma \widehat{X}_\epsilon(y) -\gamma^2 c_\epsilon}\,W(dy),$$
with
$$\Phi_\ell(x)=\frac{\varphi(x+\ell/2)}{|x+\ell/2|^{\frac{1}{2}-H}}-\frac{\varphi(x-\ell/2)}{|x-\ell/2|^{\frac{1}{2}-H}}.$$
Similarly to the computation of the variance, conditionally on the field $\widehat{X}_\epsilon$, we can use Wick's formula with respect to the white noise to get
 \begin{align*}
 \E\big[\delta_\ell u_\epsilon(x)^2|\widehat{X}_\epsilon\big] =& \int \Phi_\ell(x-y)^2\E[X_\epsilon(y)^2|\widehat{X}_\epsilon]M_{2\gamma}^\epsilon(dy)\\
 &+\iint \Phi_\ell(x-z)  \Phi_\ell(x-y)k_\epsilon(y-z)k_\epsilon(z-y)M_{2\gamma}^\epsilon(dz)M_{2\gamma}^\epsilon(dy)\\
 =&\iint \big(\Phi_\ell(x-y)^2-\Phi_\ell(x-z)  \Phi_\ell(x-y)\big)k^2_\epsilon(y-z)M_{2\gamma}^\epsilon(dz)M_{2\gamma}^\epsilon(dy).
 \end{align*}
This latter integral can be symmetrized to obtain 
\begin{equation}\label{eq:VarIncrMultibis}
\E (\delta_\ell u)^2=\lim_{\epsilon \to 0}\E (\delta_\ell u_\epsilon)^2=  \frac{1}{2} \iint  \big(\Phi_\ell(y)-\Phi_\ell(z)\big)^2  \frac{1}{|y-z|^{1+4 \gamma^2}} 1_{|y-z| \leq 1}  dz dy.
\end{equation}
This leads to 
\begin{equation}\label{eq:EquivIncrMultibis}
\E (\delta_\ell u)^2\build{\sim}_{\ell\to 0}^{} C_2 \ell^{2H-4 \gamma^2},
\end{equation}
with $C_2$ a strictly positive constant, independent of the regularization procedure, and given by
\begin{align}\label{eq:formulaC2}
C_2=  \frac{\varphi(0)^2}{2} \iint  \Big ( \frac{1}{|y+1/2|^{1/2-H}}-\frac{1}{|y-1/2|^{1/2-H}}&- \frac{1}{|z+1/2|^{1/2-H}}+ \frac{1}{|z-1/2|^{1/2-H}}  \Big )^2 \notag \\
&\times \frac{1}{|y-z|^{1+4 \gamma^2}}  dz dy.
\end{align}
We invite the reader to Section \ref{Sec:ProofIncr} devoted to the proofs of existence of the limiting value of the increment variance \eqref{eq:VarIncrMultibis} and of its equivalent at small scales \eqref{eq:EquivIncrMultibis}. Similarly to the variance, the equivalent of the increment variance \eqref{eq:EquivIncrMultibis} will eventually makes sense for $\gamma^2 < \frac 1 2 H$. As we did for the variance, we would like to present now a more straightforward way to derive the equivalent at small scales. This will also prepare for the computations of Section \ref{Sec:Skewness} below.

We have in a similar fashion the following limit of the increment variance for $\epsilon \to 0$\begin{equation}\label{eq:VarIncrMulti}
\E (\delta_\ell u)^2=\lim_{\epsilon \to 0}\E (\delta_\ell u_\epsilon)^2 =\frac{1}{2\gamma^2}\int_0^{1} (\Phi_\ell \star\Phi_\ell)'(h)\big(1-|h|^{-4\gamma^2}\big) \, dh,
\end{equation}
where,  for $H\in ]0,1[/\{1/2\}$, the derivative of the bounded function  $(\Phi_\ell \star\Phi_\ell)(h)$ is given by the following principal value integral
$$ (\Phi_\ell \star\Phi_\ell)'(h) = \mathrm{ P.V. } \int \Phi_\ell(x)\Phi_\ell'(x+h)dx.$$
As encountered in the calculation of the variance, one has to check the integrability in the neighboorhood of zero. To do so, remark that
$$ (\Phi_\ell \star\Phi_\ell)(h) =   2(\phi\star\phi)(h)-(\phi\star\phi)(h+\ell)-(\phi\star\phi)(h-\ell),$$
which leads to
$$ (\Phi_\ell \star\Phi_\ell)'(h) =   2(\phi\star\phi)'(h)-(\phi\star\phi)'(h+\ell)-(\phi\star\phi)'(h-\ell).$$
As we have seen (Lemma \ref{lem:Equiv2}), the function $(\phi\star\phi)'(h)$ diverges at the origin as fast as $\frac{1}{|h|^{1-2H}}$ when $H<1/2$ (it remains bounded when $H>1/2$), which is itself integrable in the neighboorhood of the origin. Thus the asymptotical variance of increments \eqref{eq:VarIncrMulti} makes sense as soon as $1-2H+4 \gamma^2<1$. 

Let us now compute the asymptotical behavior of the variance of increments \eqref{eq:VarIncrMulti} at vanishing scale $\ell \to 0$. We can always write for $\gamma\ne 0$
\begin{align}
\E (\delta_\ell u)^2=\frac{1}{2 \gamma^2} \int_0^{1/\ell} (\Phi_\ell \star\Phi_\ell)'(\ell h)  \big(1-|\ell h|^{-4\gamma^2}\big) \, \ell dh. \nonumber
\end{align}
Remark that
$$ (\Phi_\ell \star\Phi_\ell)(\ell h) =   2(\phi\star\phi)(\ell h)-(\phi\star\phi)(\ell (h+1))-(\phi\star\phi)(\ell (h-1)),$$
so that, using Lemma \ref{lem:Equiv2}, we have the following equivalent, 
\begin{align*}
(\Phi_\ell \star\Phi_\ell)'(\ell h)  \build{\sim}_{\ell\to 0^+}^{}  & (H-1/2) \varphi^2(0)\ell^{2H-1}  \\
 \hskip - 4 cm & \times \left[ 2\sign(h)|h|^{2H-1}-\sign(h+1)|h+1|^{2H-1}-\sign(h-1)|h-1|^{2H-1}\right]\\
 \hskip - 4 cm  &\times\mathrm{ P.V. }\int\frac{1}{|x|^{\frac{1}{2}-H}}\frac{x+1}{|x+1|^{\frac{5}{2}-H}} dx.
\end{align*}
Thus, at small scales ($\ell\to 0^+$),  for $H\in ]0,1[/\{1/2\}$, the variance of the increments behaves as
\begin{equation}\label{eq:EquivIncrMulti}
\E (\delta_\ell u)^2\build{\sim}_{\ell\to 0}^{} -\frac{a_{\gamma,H}}{2 \gamma^2} \varphi^2(0)\ell^{2H-4 \gamma^2}(H-1/2)
\times\mathrm{ P.V. }\int\frac{1}{|x|^{\frac{1}{2}-H}}\frac{x+1}{|x+1|^{\frac{5}{2}-H}} dx,
\end{equation}
where
\begin{equation*}
a_{\gamma,H}= \int_0^{\infty} \frac{1}{h^{4 \gamma^2}}\left[ 2h^{2H-1}-(h+1)^{2H-1}-\sign(h-1)|h-1|^{2H-1}\right]dh.
\end{equation*}
Remark that 
$$ \mathrm{ P.V. }\int\frac{1}{|x|^{\frac{1}{2}-H}}\frac{x+1}{|x+1|^{\frac{5}{2}-H}} dx = \int_0^{\infty}\frac{1}{x^{\frac{3}{2}-H}}\left[\frac{1}{|x-1|^{\frac{1}{2}-H}}-\frac{1}{|x+1|^{\frac{1}{2}-H}}\right]dx$$
is negative for $H>1/2$, and positive for $H<1/2$, which makes the former equivalent of $\E (\delta_\ell u)^2$ (Eq. \eqref{eq:EquivIncrMulti}) of the same sign as $a_{\gamma,H}$. 

\vspace{0.2 cm}

This approach bears a difficulty though, i.e. it is non obvious to show that $a_{\gamma,H}>0$. Fortunately, the derivation of the equivalent than we obtained with the first method \eqref{eq:EquivIncrMultibis} ensures that the multiplicative constant entering in the equivalent at small scales \eqref{eq:formulaC2} is indeed positive.

\begin{rem}
In the case $\gamma=0$, by similar computations we get the analogue of \eqref{eq:VarIncrMulti}
\begin{equation}\label{eq:VarIncrMono}
\E (\delta_\ell u)^2= 2 \int_0^{1} (\Phi_\ell \star\Phi_\ell)'(h) \ln h \, dh.
\end{equation}
One can also obtain \eqref{eq:VarIncrMono} as limit of \eqref{eq:VarIncrMulti} when $\gamma$ goes to $0$. This leads to the following equivalent
\begin{equation}\label{eq:EquivIncrMono}
\E (\delta_\ell u)^2\build{\sim}_{\ell\to 0}^{} a_{H} (H-1/2) \varphi^2(0)\ell^{2H}  \ln \ell
\times\mathrm{ P.V. }\int\frac{1}{|x|^{\frac{1}{2}-H}}\frac{x+1}{|x+1|^{\frac{5}{2}-H}} dx
\end{equation}
 where the remaining constant can be made explicit, i.e.
\begin{equation*}
a_{H}= \int_0^{\infty} \left[ 2h^{2H-1}-(h+1)^{2H-1}-\sign(h-1)|h-1|^{2H-1}\right]dh= \frac{1}{H}.
\end{equation*}
Remark that contrary to the $\gamma\ne 0$ case, in which the second order structure function $\E (\delta_\ell u)^2$ \eqref{eq:EquivIncrMulti} behaves at small scales as a power-law $\ell^{2H-4\gamma^2}$, an additional logarithmic correction appears in front of the power-law $\ell^{2H}$ in the $\gamma=0$ case \eqref{eq:EquivIncrMono}.
\end{rem}

\subsection{Skewness of increments}\label{Sec:Skewness}

For the case $\gamma>0$, using the fact that $\Phi_\ell$ is a odd function of its argument, we have by statistical homogeneity
\begin{align*}
\E  [ (\delta_\ell u_\epsilon)^3   |  \widehat{X}_\epsilon ]   = - \iiint  &\Phi_\ell(x) \Phi_\ell(y) \Phi_\ell(z)  \E[ X_\epsilon(x)X_\epsilon(y)X_\epsilon(z)   W(dx) W(dy) W(dz)  | \widehat{X}_\epsilon ]  \\
&\times e^{\gamma (\widehat{X}_\epsilon(x)+\widehat{X}_\epsilon(y)+\widehat{X}_\epsilon(z))  -3 \gamma^2 c_\epsilon}\,    dx dy dz.
\end{align*}
By standard integration by parts and exploiting symmetry $x \leftrightarrow y \leftrightarrow z $, we get (recall that $\E[ X_\epsilon(z) W(dx) | \widehat{X}_\epsilon  ]= e^{\gamma \widehat{X}_\epsilon(x)  - \gamma^2 c_\epsilon} k_\epsilon(z-x) dx$),
\begin{align*}
& \E  [ (\delta_\ell u_\epsilon)^3   |  \widehat{X}_\epsilon ]    \\
& = -6 \iint \Phi_\ell(x) \Phi_\ell(y)^2   \E[ X_\epsilon(x)X_\epsilon(y)  | \widehat{X} _\epsilon] k_\epsilon(y-x)   e^{2 \gamma (\widehat{X}_\epsilon(x)+\widehat{X}_\epsilon(y))  -4 \gamma^2 c_\epsilon}\,    dx dy   \\
& -2 \iiint \Phi_\ell(x) \Phi_\ell(y)  \Phi_\ell(z)  k_\epsilon(x-y) k_\epsilon(y-z)  k_\epsilon(z-x)     e^{2 \gamma (\widehat{X}_\epsilon(x)+\widehat{X}_\epsilon(y)+\widehat{X}_\epsilon(z) )  -6 \gamma^2 c_\epsilon}\,    dx dy dz   \\
& = -6  \iint \Phi_\ell(x) \Phi_\ell(y)^2   \E[ X_\epsilon(x)X_\epsilon(y)  | \widehat{X}_\epsilon ] k_\epsilon(y-x)   e^{2 \gamma (\widehat{X}_\epsilon(x)+\widehat{X}_\epsilon(y))  -4 \gamma^2 c_\epsilon}\,    dx dy   \\
& = -6 \iiint \Phi_\ell(x) \Phi_\ell(y)^2   k_\epsilon(x-t)   k_\epsilon(y-t)   k_\epsilon(y-x)   e^{2 \gamma (\widehat{X}_\epsilon(x)+\widehat{X}_\epsilon(y)+\widehat{X}_\epsilon(t) )  -6 \gamma^2 c_\epsilon}\,    dx dy dt,   
\end{align*}
where we have used the fact that the triple integral $ \int \cdots dxdydz$ in the above computation is equal to $0$ by symmetry. Recall that we have noted the even function $\widehat{C}_\epsilon(x)=\E[  \widehat{X}_\epsilon(x)\widehat{X}_\epsilon(0) ]$ (see \eqref{eq:DefCHat}). Now, by averaging with respect to $\widehat{X}_\epsilon$, we get 
\begin{align*}
& \E  [ (\delta_\ell u_\epsilon)^3   ]    \\
& = -6\iiint \Phi_\ell(x) \Phi_\ell(y)^2   k_\epsilon(x-t)   k_\epsilon(y-t)   k_\epsilon(y-x)  e^{4\gamma^2\left[ \widehat{C}_\epsilon(x-t)+\widehat{C}_\epsilon(y-t)+\widehat{C}_\epsilon(y-x)\right]} \,    dx dy dt  \\
& = -6 \iint  \Phi_\ell(x) \Phi_\ell(y)^2  k_\epsilon(y-x)  e^{4\gamma^2 \widehat{C}_\epsilon(y-x)} C_{\epsilon,\gamma}(y-x)  \,    dx dy   \\
& = -6 \iint  \Phi_\ell(x) \Phi_\ell(x+h)^2 k_\epsilon(h)  e^{4\gamma^2 \widehat{C}_\epsilon(h)}   C_{\epsilon,\gamma}(h)  \,    dx dh  \\
& = -6 \int   ( \Phi_\ell \star \Phi_\ell^2) (h)   k_\epsilon(h)  e^{4\gamma^2 \widehat{C}_\epsilon(h)}  C_{\epsilon,\gamma}(h)  \,     dh,
\end{align*}
where we defined analogously the even function $C_{\epsilon,\gamma}(h) $ explicitly given by
\begin{equation*} 
C_{\epsilon,\gamma}(h)= \int k_\epsilon(x) k_\epsilon(x+h)e^{4\gamma^2\left[ \widehat{C}_\epsilon(x)+\widehat{C}_\epsilon(x+h)\right]}   dx.
\end{equation*}
We have pointwise for $h\ne 0$
\begin{equation*} 
C_\gamma (h):= \lim_{\epsilon\to 0}C_{\epsilon,\gamma}(h) =
\begin{cases}
& 2\ln_+\left(\frac{1}{|h|}\right)+g_0(h), \:  \text{if} \: \gamma=0 \\
&  \frac{r_\gamma}{|h|^{8\gamma^2}}+g_\gamma(h), \:   \text{if} \: \gamma^2\in]0,1/8[\\
\end{cases}
\end{equation*}
with $g_\gamma$ a bounded function of its argument for any $\gamma\ge 0$, and we have set for $\gamma^2\in]0,1/8[$  
\begin{equation*}
r_\gamma= 2\int_0^\infty  \frac{1}{\sqrt{x}|x|^{4\gamma^2} } \Big(\frac{x+1}{|x+1|^{\frac{3}{2}+4\gamma^2}}+\frac{x-1}{|x-1|^{\frac{3}{2}+4\gamma^2}}\Big)dx.
\end{equation*}

Hence, we get the identity
\begin{equation}\label{eq:Mom3IncrMulti}
 \E  [ (\delta_\ell u)^3   ] := \underset{\epsilon \to 0}{\lim}  \:  \E  [ (\delta_\ell u_\epsilon)^3   ]  =  -12 \int_0^1   (\Phi_\ell \star \Phi_\ell^2) (h)   \frac{  1}{h^{\frac{1}{2}+4 \gamma^2}}   C_{\gamma}(h)  \,     dh.
\end{equation}

To make sense of the of the asymptotical form of the third moment of increments \eqref{eq:Mom3IncrMulti}, we have to check the integrability of the integrand in the neighborhood of the origin. To do so, we have to study the behavior of the function $ (\Phi_\ell \star \Phi_\ell^2)$: this is the subject of Lemma \ref{lem:PhiPhi2Bound}. We show there that the function is singular at $h=\ell$ only in the case $H\in ]0,1/6]$,  a singular behavior that is itself integrable. For $H\in ]1/6,1[/\{1/2\}$, $(\Phi_\ell \star \Phi_\ell^2)$ is a continuous and bounded function of its argument. At the origin, $\forall H$, $(\Phi_\ell \star \Phi_\ell^2)(h)$ goes to zero as fast as $h$. Thus, the equivalent \eqref{eq:Mom3IncrMulti} makes sense for $H\in ]0,1[/\{1/2\}$ and $\gamma^2<1/8$.

\begin{lemma}\label{lem:PhiPhi2Bound}
For $H\in ]1/6,1[$ and $\forall h$, $(\Phi_\ell\star\Phi^2_\ell)(h)$ is a continuous and bounded function of its argument. For $H\in ]0,1/6]$, $(\Phi_\ell\star\Phi^2_\ell)$ has an additional singularity at $h=\ell$ given by
\begin{equation*}
(\Phi_\ell\star\Phi^2_\ell)(h) \underset{h \to \ell}{\sim}
\begin{cases}
& d_H\varphi^3(0) |h-\ell|^{3H-\frac{1}{2}}, \;  \text{if} \;  H<1/6\\
& 2\varphi^3(0) \ln \frac{1}{|h-\ell|}, \;   \text{if} \; H=1/6, \\
\end{cases}
\end{equation*}
where $d_H$ is a constant independent of $\varphi(0)$ that we can compute.
Furthermore, for any $H\in ]0,1[/\{1/2\}$ we have the following equivalent at small arguments
\begin{align*}
(\Phi_\ell\star\Phi^2_\ell)(h)\build{\sim}_{h\to 0}^{} h&\int \left[ (1/2-H)\varphi(x)-x\varphi'(x)\right]\frac{x}{|x|^{5/2-H}}\left[\frac{\varphi(x-\ell)}{|x-\ell|^{1/2-H}}-\frac{\varphi(x+\ell)}{|x+\ell|^{1/2-H}}\right]\\
&\times \left[\frac{\varphi(x-\ell)}{|x-\ell|^{1/2-H}}+\frac{\varphi(x+\ell)}{|x+\ell|^{1/2-H}}-\frac{2\varphi(x)}{|x|^{1/2-H}}\right]dx.
\end{align*} 

\end{lemma}
 \noindent {\it Proof.} We have
\begin{align*} 
(\Phi_\ell\star\Phi^2_\ell)(h) = \int&\left[\frac{\varphi(x+\ell/2)}{|x+\ell/2|^{1/2-H}}-\frac{\varphi(x-\ell/2)}{|x-\ell/2|^{1/2-H}}\right]\\
&\times\left[\frac{\varphi(x+h+\ell/2)}{|x+h+\ell/2|^{1/2-H}}-\frac{\varphi(x+h-\ell/2)}{|x+h-\ell/2|^{1/2-H}}\right]^2dx.
\end{align*} 
 Notice that $(\Phi_\ell\star\Phi^2_\ell)(h)$ can be written with the following convenient form
\begin{align*} 
(\Phi_\ell\star\Phi^2_\ell)(h) = \int\frac{\varphi(x-h)}{|x-h|^{1/2-H}}&\left[\frac{\varphi(x-\ell)}{|x-\ell|^{1/2-H}}-\frac{\varphi(x+\ell)}{|x+\ell|^{1/2-H}}\right]\\
&\times \left[\frac{\varphi(x-\ell)}{|x-\ell|^{1/2-H}}+\frac{\varphi(x+\ell)}{|x+\ell|^{1/2-H}}-\frac{2\varphi(x)}{|x|^{1/2-H}}\right]dx,
\end{align*} 
which shows that $(\Phi_\ell\star\Phi^2_\ell)$ is continuous and bounded for $H\in ]1/6,1[$ and $\forall h$. For $H\in ]0,1/6]$,  $(\Phi_\ell\star\Phi^2_\ell)$ has an additional singularity at $h=\ell$. The proposed equivalent for $h\to 0$ follows from the Taylor Series of the first ratio entering in the integral (the first contribution to this development vanishes by symmetry).

Let us now take a look at the additional singularity when $H\le 1/6$. From this former expression, we see that $(\Phi_\ell\star\Phi^2_\ell)(h)$ as the same singularity when $h$ goes to $\ell$ as
\begin{equation*}
\int_{|x| \leq 1}  \frac{\varphi(x+\ell-h)}{|x+\ell-h|^{1/2-H}}\frac{\varphi^2(x)}{|x|^{1-2H}}dx.
\end{equation*}
If $H< \frac{1}{6}$ then it is equal to, take for instance $h<\ell$,
\begin{align*}
|h-\ell|^{3H-1/2}\int_{|y| \leq \frac{1}{|h-\ell|}} &\frac{\varphi[|h-\ell|(y+1)]}{|y+1|^{1/2-H}}\frac{\varphi^2[|h-\ell|y]}{|y|^{1-2H}}dy\\
&\underset{|h-\ell| \to 0}{\sim} d_H\varphi^3(0)|h-\ell|^{3H-\frac{1}{2}},
\end{align*}
where 
$$d_H=\int_{\mathbb{R}}  \frac{1}{|y+1|^{\frac{1}{2}-H }|y|^{1-2H}}  dy.$$
If $H= \frac{1}{6}$ then it is equal to
\begin{align*}
\int_{|y| \leq \frac{1}{|h-\ell|}} \frac{\varphi[|h-\ell|(y+1)]}{|y+1|^{1/3}}\frac{\varphi^2[|h-\ell|y]}{|y|^{2/3}}dy\underset{|h-\ell| \to 0}{\sim} 2 \varphi^3(0)\ln \frac{1}{|h-\ell|}.
\end{align*}
 \qed

This shows that the asymptotical form \eqref{eq:Mom3IncrMulti} makes sense for $H\in ]0,1[/\{1/2\}$ and $\gamma^2<1/8$. Let us now compute its behavior in the limit of vanishing scales $\ell \to 0$. We can always write
 \begin{align*}
\E (\delta_\ell u)^3&=-12\int_0^{1/\ell} (\Phi_\ell\star\Phi^2_\ell)(\ell h)\left( \frac{r_\gamma}{| \ell h  |^{8 \gamma^2}} +g_\gamma(\ell h) \right)\frac{1}{(\ell h)^{\frac{1}{2}+4 \gamma^2}} \ell dh,
\end{align*}
so if the following integral makes sense, we obtain the equivalent at small scales
\begin{align}\label{eq:EquivMom3SmallScales}
\E (\delta_\ell u)^3&\build{\sim}_{\ell \to 0}^{} -12 r_\gamma \varphi^3(0)\ell^{3H-12\gamma^2} \int_0^{\infty} f_H(h)\frac{1}{h^{\frac{1}{2}  +12 \gamma^2 }}dh,
\end{align}
where 
\begin{align}\label{deff_H}
f_H(h)= \int \left[\frac{1}{|x+1/2|^{\frac{1}{2}-H}}-\frac{1}{|x-1/2|^{\frac{1}{2}-H}}\right]&\left[\frac{1}{|x+h+1/2|^{\frac{1}{2}-H}}-\frac{1}{|x+h-1/2|^{\frac{1}{2}-H}}\right]^2dx.
\end{align}
To make sense of the equivalent we wrote in \eqref{eq:EquivMom3SmallScales}, in a similar manner as we did for the asymptotical form of the third moment of increments \eqref{eq:Mom3IncrMulti}, we have to check the integrability of the proposed integrand. This is the subject of Lemma \ref{lem:Equiv3}. We show there that similarly the function is singular at $h=1$ only in the case $H\in ]0,1/6]$. For $H\in ]1/6,1[/\{1/2\}$, $f_H$ is a continuous and bounded function of its argument. As far as integrability at large $h$ is concerned (we lost in this limit the cut-off function $\varphi)$, we show that for $H\in ]0,1[/\{1/2\}$, $f_H(h)$ decreases as fast as $1/h^{3/2-H}$ which is integrable when weighted by the factor $1/h^{1/2+12\gamma^2}$, for any $\gamma\ge 0$. At the origin, once again, $\forall H$, $f_H(h)$ goes to zero as fast as $h$. Thus, the equivalent \eqref{eq:EquivMom3SmallScales} makes sense also for $H\in ]0,1[/\{1/2\}$ and $\gamma^2<1/8$.

 \begin{lemma}\label{lem:Equiv3}
For $H\in ]1/6,1[/\{1/2\}$ and $\forall h$, $f_H$ is a continuous and bounded function of its argument. For $H\in ]0,1/6]$, $f_H$ has an additional singularity at $h=1$ given by
\begin{equation*}
f_H(h) \underset{h \to 1}{\sim}
\begin{cases}
& d_H |h-1|^{3H-\frac{1}{2}}, \;  \text{if} \;  H<1/6\\
& 2 \ln \frac{1}{|h-1|}, \;   \text{if} \; H=1/6, \\
\end{cases}
\end{equation*}
where $d_H$ is the same constant entering in Lemma \ref{lem:PhiPhi2Bound}. Furthermore, for $H\in ]0,1[/\{1/2\}$ we have the following equivalent at small arguments
\begin{align*}
f_H(h)\build{\sim}_{h\to 0}^{} -(H-1/2)h&\int \frac{x}{|x|^{5/2-H}}\left[\frac{1}{|x-1|^{1/2-H}}-\frac{1}{|x+1|^{1/2-H}}\right]\\
&\times \left[\frac{1}{|x-1|^{1/2-H}}+\frac{1}{|x+1|^{1/2-H}}-\frac{2}{|x|^{1/2-H}}\right]dx,
\end{align*} 
and the following equivalent at large arguments
\begin{align*}
f_H(h)\build{\sim}_{h\to \infty}^{} -(H-1/2)h^{H-3/2}&\int x\left[\frac{1}{|x-1|^{1/2-H}}-\frac{1}{|x+1|^{1/2-H}}\right]\\
&\times \left[\frac{1}{|x-1|^{1/2-H}}+\frac{1}{|x+1|^{1/2-H}}-\frac{2}{|x|^{1/2-H}}\right]dx.
\end{align*} 
 \end{lemma}

 \noindent {\it Proof.} Noticing once again that $f_H$ can be written with the following convenient form
\begin{align*} 
f_H(h) = \int\frac{1}{|x-h|^{1/2-H}}&\left[\frac{1}{|x-1|^{1/2-H}}-\frac{1}{|x+1|^{1/2-H}}\right]\\
&\times \left[\frac{1}{|x-1|^{1/2-H}}+\frac{1}{|x+1|^{1/2-H}}-\frac{2}{|x|^{1/2-H}}\right]dx,
\end{align*} 
proofs are then similar to those of Lemma \ref{lem:PhiPhi2Bound}. The proposed equivalent $h\to \infty$ follows from the factorization of $h$ in the first ratio and then doing a Taylor Series.
 \qed

\begin{rem}
It remains to show that indeed the integral entering in the equivalent \eqref{eq:EquivMom3SmallScales} does not vanish: this is a direct consequence of the assumption \eqref{conjf_H}. Indeed, as illustrated by a numerical estimation (see Annex \ref{ann:NumEstF_H}), the function $f_H(h)$ seems to be, $\forall h>0$, strictly positive for $H<1/2$, and strictly negative for $H>1/2$, which makes the equivalent \eqref{eq:EquivMom3SmallScales} non vanishing.
\end{rem}

\begin{rem}
Concerning the modeling of fluid turbulence, let us take a look at the predictions of the present stochastic model. We recall that empirical estimations give $4\gamma^2=0.025$ \cite{Fri95,CheCas12} and that a statistical property of stationary solutions of forced Navier-Stokes equations, namely the  $4/5^{\mbox{\tiny th}}$-law of turbulence \cite{Fri95} (see also Section \ref{Physicsmotivation}) gives
$$ \E (\delta_\ell u)^3 \build{\sim}_{\ell\to 0}^{} -\frac{4}{5}\varepsilon \ell,$$
where $\varepsilon$ is the average viscous dissipation per unit of mass (see  \cite{Fri95} for a precise definition). Taking $H=1/3+4\gamma^2$, we see here that the present model indeed predicts that $\varepsilon$ becomes independent of the viscosity, as required by the second axiom of Kolmogorov' phenomenology depicted in Section \ref{Physicsmotivation}. To see this analogy, assume that the scale $\epsilon$ entering in the regularization of the field $u_\epsilon$ plays the role of the dissipative  length scale, that is expected to go to zero as viscosity goes to 0. Furthermore, the model, as it is defined, gives the correct sign for the third order structure function \eqref{eq:EquivMom3SmallScales}, if we assume that for this $H$, $f_H(h)>0$ for $h>0$, as it is assumed in \eqref{conjf_H}, and confirmed in Annex \ref{ann:NumEstF_H}.
\end{rem}

\begin{rem}
In the case $\gamma=0$, the identity \eqref{eq:Mom3IncrMulti} is valid by simply setting $\gamma=0$. This leads to the following equivalent by similar computations to the $\gamma>0$ case
\begin{align}\label{eq:EquivMom3SmallScalesGamma0}
\E (\delta_\ell u_\epsilon)^3&\build{\sim}_{\ell \to 0}^{}-24 \varphi^3(0)\ell^{3H}\ln\left( \frac{1}{\ell}\right)\int_0^{\infty} f_H(h)\frac{1}{\sqrt{h}}dh.
\end{align}
\end{rem}

\subsection{High order moments}\label{ss.mom}
  
 Notice that for $\epsilon>0$, and conditionally on the field $\widehat{X}$, the field $u_\epsilon$ is in the second Wiener chaos generated by the white noise $W$. Hence by  hyper-contractivity in the second Wiener chaos  (see Theorem 2.7.2 and Corollary 2.8.14  in \cite{NouPec12}), we get the existence of some constant $C_q>0$ such that, for $q\geq 0$,
 $$C_q^{-1}\E\Big[u_\epsilon(x)^2|\widehat{X}_\epsilon\Big]^{q/2}\leq \E\Big[|u_\epsilon(x)|^q|\widehat{X}_\epsilon\Big]\leq C_q\E\Big[u_\epsilon(x)^2|\widehat{X}_\epsilon\Big]^{q/2}.$$
 Therefore
 $$C_q^{-1} \E\Big[\E\big[u_\epsilon(x)^2|\widehat{X}_\epsilon\big]^{q/2}\Big]\leq \E[|u_\epsilon(x)|^q]\leq C_q\E\Big[\E\big[u_\epsilon(x)^2|\widehat{X}_\epsilon\big]^{q/2}\Big].$$
Recall that we have  
\begin{align*}
 \E\big[u_\epsilon(x)^2|\widehat{X}_\epsilon\big] =&\frac{1}{2}\iint \big(\phi(x-y)-\phi(x-u)\big)^2k^2_\epsilon(y-u)M_{2\gamma}^\epsilon(du)M_{2\gamma}^\epsilon(dy).
 \end{align*}
Now use the Gaussian multiplicative chaos technology to check the integrability properties of this integral.  By statistical homogeneity, it suffices to consider the case $x=0$ and, as a result of the above discussion, it is enough to show that
 $$\sup_{\epsilon\in]0,1]}\E\Big[\Big(\iint \big(\phi(y)-\phi(u)\big)^2k^2_\epsilon(y-u)M_{2\gamma}^\epsilon(du)M_{2\gamma}^\epsilon(dy)\Big)^{q/2} \Big]<+\infty.$$
 From now on we assume   $q\geq2$ and we will not treat the case $q<2$ as it is quite similar: we will only mention below how to adapt the proof.
 The latter quantity is finite provided that we can show that $\sup_{\epsilon\in]0,1]} A_\epsilon([0,1]^2)<+\infty$ with 
 $$A_\epsilon(D):=\E\Big[\Big(\iint_{D} \big(\frac{1}{|y|^{\tfrac{1}{2}-H}}-\frac{1}{|u|^{\tfrac{1}{2}-H}}\big)^2\frac{1}{|u-y|\vee \epsilon}M_{2\gamma}^\epsilon(du)M_{2\gamma}^\epsilon(dy)\Big)^{q/2} \Big]^{2/q}.$$
By Kahane's convexity inequality \cite{Kah85} (see also \cite{RhoVar14}) we may and will assume that $\widehat{X}_\epsilon$ is the exact scale invariant kernel studied in \cite{BacMuz03}.
Hence, for all $\lambda,\epsilon \in]0,1[$, it satisfies the following equality in law  
\begin{equation}\label{exact}
(\widehat{X}_{\lambda\epsilon }(\lambda u))_{u\in [-1,1]}=(\widehat{X}_{\epsilon}(u)+\Omega_\lambda)_{u\in [-1,1]}
\end{equation}
 where $\Omega_\lambda$ is a centered Gaussian random variable with variance $-\ln \lambda$   independent of the process $\widehat{X}_{\epsilon}$. We have
\begin{align*}
A_{\epsilon/2}([0,1]^2)&\leq A_{\epsilon/2}([0,1/2]^2)\\
&+A_{\epsilon/2}([0,1/2]\times[1/2,1])+A_{\epsilon/2}([1/2,1]\times[0,1/2])+A_{\epsilon/2}([1/2,1]^2).
\end{align*} 
Here we mention that we have used the triangular inequality for $L^p$ norms with $p=q/2\geq 1$ (in the case $q<2$, just use the subadditivity of the mapping $x\mapsto x^{q/2}$).
Thanks to Corollary \ref{lem2} below, we deduce that, for some irrelevant constant $C>0$ and all $\epsilon\in]0,1]$
\begin{equation}\label{recA}
A_{\epsilon/2}([0,1]^2)\leq A_{\epsilon/2}([0,1/2]^2)+C.
\end{equation}
Let us make the changes of variables $u'=2u$ and $y'=2y$, we get
 \begin{align*}
 A_{\epsilon/2}([0,1/2]^2)=2^{-2H}\E\Big[I_\epsilon^{q/2} \Big]^{2/q},
 \end{align*}
where
$$I_\epsilon= \iint_{[0,1]^2} \big(\frac{1}{|y'|^{\tfrac{1}{2}-H}}-\frac{1}{|u'|^{\tfrac{1}{2}-H}}\big)^2\frac{1}{|u'-y'|\wedge\epsilon}e^{2\gamma \widehat{X}_{\epsilon/2}(u'/2)+2\gamma \widehat{X}_{\epsilon/2}(y'/2)-4\gamma^2\E[\widehat{X}_{\epsilon/2}^2]} du'dy'.$$
 Now we can use the relation in law \eqref{exact} to get
 \begin{align*}
 A_{\epsilon/2}([0,1/2]^2)&=2^{-2H}\E\big[e^{2\gamma q\Omega_{1/2}}\big]^{2/q}e^{-4\gamma^2\ln 2}\times\\
 &\E\Big[\Big(\iint_{[0,1]^2} \big(\frac{1}{|y'|^{\tfrac{1}{2}-H}}-\frac{1}{|u'|^{\tfrac{1}{2}-H}}\big)^2\frac{1}{|u'-y'|\wedge\epsilon}M_{2\gamma}^\epsilon(du')M_{2\gamma}^\epsilon(dy')\Big)^{q/2} \Big]^{2/q}\\
 =&2^{-2H-4\gamma^2+4\gamma^2 q}  A_{\epsilon}([0,1]^2).
 \end{align*}
 Under the assumption $q<1+\tfrac{H}{2\gamma^2}$, the exponent of $2$ in the above expression is strictly negative. From \eqref{recA}, we deduce
$$A_{\epsilon/2}([0,1]^2)\leq r A_{\epsilon}([0,1]^2)+C$$ for some constant $r\in]0,1[$,
  ensuring finiteness of the supremum of the family $(A_{\epsilon}([0,1]^2))_\epsilon$ as claimed.\qed

 Before proving Corollary \ref{lem2}, we recall the following lemma which is a 1d version of lemma A.1 in \cite{DaKu16} (the proof follows the same argument as lemma A.1 in \cite{DaKu16}):
 
\begin{lemma} \label{lem1}
For $\alpha\in [0,1+2\gamma^2[$ and $q\in[0,\tfrac{1+2\gamma^2-\alpha}{2\gamma^2}\wedge \tfrac{1}{2\gamma^2}[$, we have
 $$\sup_{\epsilon\in]0,1]}\E\Big[\Big(\int_{[0,1]}  \frac{1}{(|y|\vee \epsilon)^{\alpha}}  M_{2\gamma}^\epsilon(dy)\Big)^{q} \Big]<+\infty$$
\end{lemma} 

With this lemma we can now prove the following Corollary:

 \begin{corollary} \label{lem2}
For  $H\in ]0,1/2[$  and $q\in[0,(1+\tfrac{H}{\gamma^2})\wedge \tfrac{1}{2 \gamma^2}[$, we have
\begin{align*}
1)&\sup_{\epsilon\in]0,1]}\E\Big[\Big(\iint_{[0,1/2]\times [1/2,1]} \big(\frac{1}{|y|^{\tfrac{1}{2}-H}}-\frac{1}{|u|^{\tfrac{1}{2}-H}}\big)^2 \frac{1}{|y-u|\vee \epsilon}  M_{2\gamma}^\epsilon(dy)M_{2\gamma}^\epsilon(du)\Big)^{q/2} \Big]^{2/q}<+\infty,\\
2)&\sup_{\epsilon\in]0,1]}\E\Big[\Big(\iint_{  [1/2,1]^2} \big(\frac{1}{|y|^{\tfrac{1}{2}-H}}-\frac{1}{|u|^{\tfrac{1}{2}-H}}\big)^2 \frac{1}{|y-u|\vee \epsilon}  M_{2\gamma}^\epsilon(dy)M_{2\gamma}^\epsilon(du)\Big)^{q/2} \Big]^{2/q}<+\infty.
\end{align*}
\end{corollary} 

\noindent {\it Proof.} (fix $q>2$ otherwise use sub-additivity) We can divide the square $[0,1/2]\times [1/2,1]$ into two pieces: $[0,1/4]\times [1/2,1]$ and $[1/4,1/2]\times [1/2,1]$. 

The above supremum when integrating over $[0,1/4]\times [1/2,1]$ is obviously less than
$$C\sup_{\epsilon\in]0,1]}\E\Big[\Big(\iint_{[0,1/4]\times [1/2,1]} \big(1+\frac{1}{|u|^{1-2H}}\big)    M_{2\gamma}^\epsilon(dy)M_{2\gamma}^\epsilon(du)\Big)^{q/2} \Big]^{2/q}$$ 
for some irrelevant constant $C>0$. This quantity is again less than   (up to irrelevant multiplicative constant) 
$$\sup_{\epsilon\in]0,1]}\E\Big[M_{2\gamma}^\epsilon([1/2,1])^{q}\Big]^{2/q}+ \sup_{\epsilon\in]0,1]}\E\Big[ \Big(\int_{[0,1/4]} \big(1+\frac{1}{|u|^{1-2H}}\big)     M_{2\gamma}^\epsilon(du)\Big)^{q} \Big]^{2/q}.$$ 
Indeed, this can be shown for $q\geq 2$  by making use of the elementary inequality $ab\leq a^2/2+b^2/2$ and then Minkowski inequality and for $q\leq 2$, by sub-additivity of the mapping $x\in\R_+\mapsto x^{q/2}$ and then convexity of the mapping $x\in\R_+\mapsto x^{2/q}$.  
The first supremum is finite if $q<1/(2\gamma^2)$ by  standard results on GMC theory. The second one is finite provided that $q\in [0,  1+\tfrac{H}{\gamma^2}[$.

On the area $[1/4,1/2]\times [1/2,1]$, we have the inequality
\begin{equation*}
\big(\frac{1}{|y|^{\tfrac{1}{2}-H}}-\frac{1}{|u|^{\tfrac{1}{2}-H}}\big)^2 \frac{1}{|y-u|\vee \epsilon} \leq C |y-u|  \leq  C
\end{equation*} 
so that the corresponding supremum can be shown to be finite for $q<1/(2 \gamma^2)$ by  standard results on GMC theory again. The latter argument also holds for our second claim.\qed

\subsection{Multifractal spectrum}\label{Sec:MultiSpecEnca}
Recall that we have set 
$$\delta_\ell u_\epsilon (x) = u_\epsilon (x+\ell/2)-u_\epsilon (x-\ell/2)=\int _{\mathbb R}\Phi_\ell(x-y) X_\epsilon(y)e^{\gamma\hat{X}_\epsilon(y)-\gamma^2\E[\hat{X}_\epsilon^2]}\,W(dy),$$
with
$$\Phi_\ell(x)=\frac{\varphi(x+\ell/2)}{|x+\ell/2|^{\frac{1}{2}-H}}-\frac{\varphi(x-\ell/2)}{|x-\ell/2|^{\frac{1}{2}-H}}.$$
 For $H\in ]0,1/2[$ and $0\leq q<\tfrac{1}{2\gamma^2}\wedge (1+\tfrac{H}{2\gamma^2})$, one has 
$$\E[|\delta_\ell u_\epsilon (x)| ^q]:=\lim_{\epsilon\to 0}\E[|\delta_\ell u_\epsilon (x)| ^q]<+\infty.$$
Furthermore
$$  \E[|\delta_\ell u  (x)| ^q]\sim_q \E\Big[\Big(\iint \big(\Phi_\ell(x-y)-\Phi_\ell(x-u)\big)^2k^2 (y-u)M_{2\gamma} (du)M_{2\gamma} (dy)\Big)^{q/2}\Big]  $$
where $A\sim_q B$ means here that there exists $c_q>0$ (constant only depending on $q$) such that $c_q^{-1}A\leq B\leq c_qA$.

\begin{proposition}\label{prop:spectrum}
For $H\in ]0,1/2[$ and $0\leq q<\tfrac{1}{2\gamma^2}\wedge ( 1+\tfrac{H}{2\gamma^2})$, we have for $\ell \in]0,1]$
$$ \E[|\delta_\ell u  (x)| ^q]\sim_q \ell^{\left(H+2\gamma^2\right)q-2\gamma^2q^2}.$$
\end{proposition}

 \subsection{Continuity of the limiting process $u$}

\noindent {\it Proof of continuity of $u$ in Theorem \ref{th1secdef}.}

\ni
For $H\in ]0,1[$ and $\gamma^2<1/2$,  we have shown tightness of the finite dimensional marginals via moment estimates. The limiting process $u$ satisfies for $0\leq q<\tfrac{1}{2\gamma^2}\wedge ( 1+\tfrac{H}{2\gamma^2})$
$$\forall x,y\in\R,\quad \E[(u(x)-u(y))^q]\leq C_q|x-y|^{\left(H+2\gamma^2\right)q-2\gamma^2q^2}$$  
Furthermore, because $\gamma^2<1/2$,  setting $q_0:=\frac{1}{\sqrt{2}\gamma}$, we have $q_0<\tfrac{1}{2\gamma^2}\wedge (  1+\tfrac{H}{2\gamma^2})$. Notice that $\left(H+2\gamma^2\right)q_0-2\gamma^2q_0^2>1$ because $H+(\sqrt{2}\gamma-1)^2>1$. Hence Kolmogorov's continuity criterion ensures that $u$ admits a continuous modification such that, almost surely, its sample paths on any compact interval are $\alpha$-Holder for any $\alpha<\frac{\left(H+2\gamma^2\right)q_0-2\gamma^2q_0^2-1}{q_0}=H+ (\sqrt{2}\gamma-1)^2-1 $.\qed

\section{Analysis of the $q^{{\tiny th}}$ moments of increments} \label{Sec:ProofIncr}

This section is devoted to the rigorous derivation of the equivalents of the increments at small scales $\ell \to 0$. In Section \ref{Sec:MultiSpecEnca} (see Proposition \ref{prop:spectrum}), we have been able to show that the $q^{{\tiny th}}$-order structure function, namely $ \E[|\delta_\ell u  (x)| ^q]$, is bounded from above and below by a power-law. We would like here to go further and compute the precise equivalent.

For the sake of generality, we will perform all the calculations with an additional parameter, say $\tilde H$, that will enter in the deterministic kernel $k_\epsilon$, that we will call $k_{\epsilon,\tilde H}$. The present model $u_\epsilon$ \eqref{eq:SkMulProc} would eventually be obtained while taking $\tilde H=0$, and more generally, we will take it to be small, i.e. $\tilde H<q\gamma^2$ where the $q$ is the order of the increments moment, to ensure the intermittent corrections that are proposed in Theorem \ref{th1secdef}. 

Let us then recall the model that we will be studying here: we consider the 1d velocity field
$$u_\epsilon(x):=\int  \phi(x-y) X_\epsilon(y) e^{\gamma \widehat{X}_\epsilon(y) -\gamma^2 c_\epsilon}\,W(dy)$$
where 
\[
\begin{cases}
X_\epsilon(x)&: =\int k_{\epsilon,\tilde H}(x-y)e^{\gamma \widehat{X}_\epsilon(y) -\gamma^2 c_\epsilon}\,W(dy) \\
\widehat{C}_\epsilon(x)& := \E[\widehat{X}_\epsilon(x)\widehat{X}_\epsilon(0)]\sim  \ln_+\frac{1}{|x|+\epsilon}\\
c_\epsilon&:=\E[\widehat{X}_\epsilon(x)^2] \\
k_{\epsilon,\tilde H}(u)&:=\frac{x}{|x|_\eps^{3/2- \tilde{H}}}1_{|x|\le 1} \\
\phi(x)& :=\varphi(x)\frac{1}{|x|^{\frac{1}{2}-H}}
\end{cases}
\]
with $\varphi$ a $\mathcal C^\infty$ cut-off function of characteristic size $1$, compactly supported for instance, we choose is to be even $\varphi(x)=\varphi(-x)$ and typically one can choose for example $\varphi(x)=e^{-x^2/(1-x^2)} 1_{|x|\le 1}$.  We have introduced a new parameter $\tilde{H}\geq 0$, which has to be thought of as being small (see assumptions on $\tilde{H}$ below).

\paragraph{Summary of results with  $\tilde H$.}

This section is divided into two subsections:
\bi
\item The first one, Subsection \ref{ss.var}, which as to be seen as a warm-up, handles the case $q=2$. It has two motivations: a) introducing some spatial decomposition useful for $q\geq 2$ and b) it allows to double check the proof for $q\geq 2$. 

The conclusion of Subsection \ref{ss.var} is that under the following set of constraints:
\bnum
\item \textbf{Assumption 1}: $(q=)2 < \frac{H + \tilde H}{\gamma^2} \wedge \frac{1+\tilde H}{\gamma^2}$  \,\,\, (  $\equiv \, 2H + 2\tilde H - 4\gamma^2 >0$
 and $2\gamma^2 < 1 +\tilde H$)
\item \textbf{Assumption 2}:   $\tilde H < 2 \gamma^2$   \,\,\quad (  implies in particular $  2H + 2\tilde H - 4\gamma^2 <2$)
\enum
we have as $\ell \to 0$ 
\begin{align*}\label{}
\Eb{(\delta_\ell u)^2} \asymp  \ell^{2(H+\tilde H) - 4 \gamma^2} .
\end{align*}

\item The second one, subsection \ref{ss.NLS} deals with $q\geq 2$.  
The results of this subsection now read as follows. Under the following set of constraints:
\bnum
\item \textbf{Assumption 1}: $q< (1+ \frac {H+\tilde H} {2\gamma^2}) \wedge  \frac 1 {2\gamma^2}$

\item \textbf{Assumption 2}: $2 H + 2\tilde H - 4\gamma^2 <2$
\item \textbf{Assumption 3:}  $\tilde H < q\gamma^2 $    
\enum
we have as $\ell \to 0$,
\begin{align*}
\Eb{|\delta_\ell u|^q}   \asymp \ell^{q(H + \tilde H + 2 \gamma^2) - 2 q^2 \gamma^2 } .
\end{align*}
\ei

Recall now the expression for the velocity increments. As we have seen in Section \ref{Sec:VarVeloIncr}, we have
$$\delta_\ell u_\epsilon(x) = u_\epsilon (x+\ell/2)-u_\epsilon (x-\ell/2)=\int _{\mathbb R}\Phi_\ell(x-y)  X_\epsilon(y) e^{\gamma \widehat{X}_\epsilon(y) -\gamma^2 c_\epsilon}\,W(dy)\,,$$
with
$$\Phi_\ell(x)=\frac{\varphi(x+\ell/2)}{|x+\ell/2|^{\frac{1}{2}-H}}-\frac{\varphi(x-\ell/2)}{|x-\ell/2|^{\frac{1}{2}-H}}.$$

\ni
Following the same approach as in Sections \ref{Sec:VarVeloIncr} and  \ref{ss.mom}, we are lead to estimate:
  \begin{align*}
 \E\big[(\du_\epsilon(x))^2|\hat{X}_\epsilon\big] =&\frac{1}{2}\iint \big(\Pl(x-y)-\Pl(x-u)\big)^2k^2_\epsilon(y-u)M_{2\gamma}^\epsilon(du)M_{2\gamma}^\epsilon(dy).
 \end{align*}
From the definition of $\Pl$, it is straightforward to check that if $q\geq 2$,  
\begin{align*}
& \Eb{|\du_\eps(x=0)|^q}=\Eb{ \Eb{|\du_\eps|^q \md \hat{X}_\eps}}\\
& \leq C_q \Eb{ \Eb{(\du_\eps)^2 \md \hat{X}_\eps}^{q/2}} 
\text{ 
   (because conditioned on $\hat{X}_\eps$, we are still in the second Wiener chaos)
   }\\
& \leq \tilde C_q \Eb{
J_\epsilon^{q/2}
},
\end{align*}
where 
$$ J_\epsilon=\iint_{[-1,1]^2} 
\big(\frac{1}{|y+ \frac \ell 2|^{\tfrac{1}{2}-H}}
- \frac{1}{|y- \frac \ell 2|^{\tfrac{1}{2}-H}}
- \frac{1}{|u + \frac \ell 2|^{\tfrac{1}{2}-H}}
+ \frac{1}{|u - \frac \ell 2|^{\tfrac{1}{2}-H}}
\big)^2\frac{1}{|u-y|_\eps^{1-  2\tilde H}}M_{2\gamma}^\epsilon(du)M_{2\gamma}^\epsilon(dy).
$$

\subsection{Variance of increments $\du_\eps$}\label{ss.var}
 
\begin{figure}[h]
\begin{center}
\includegraphics[width=0.75\textwidth]{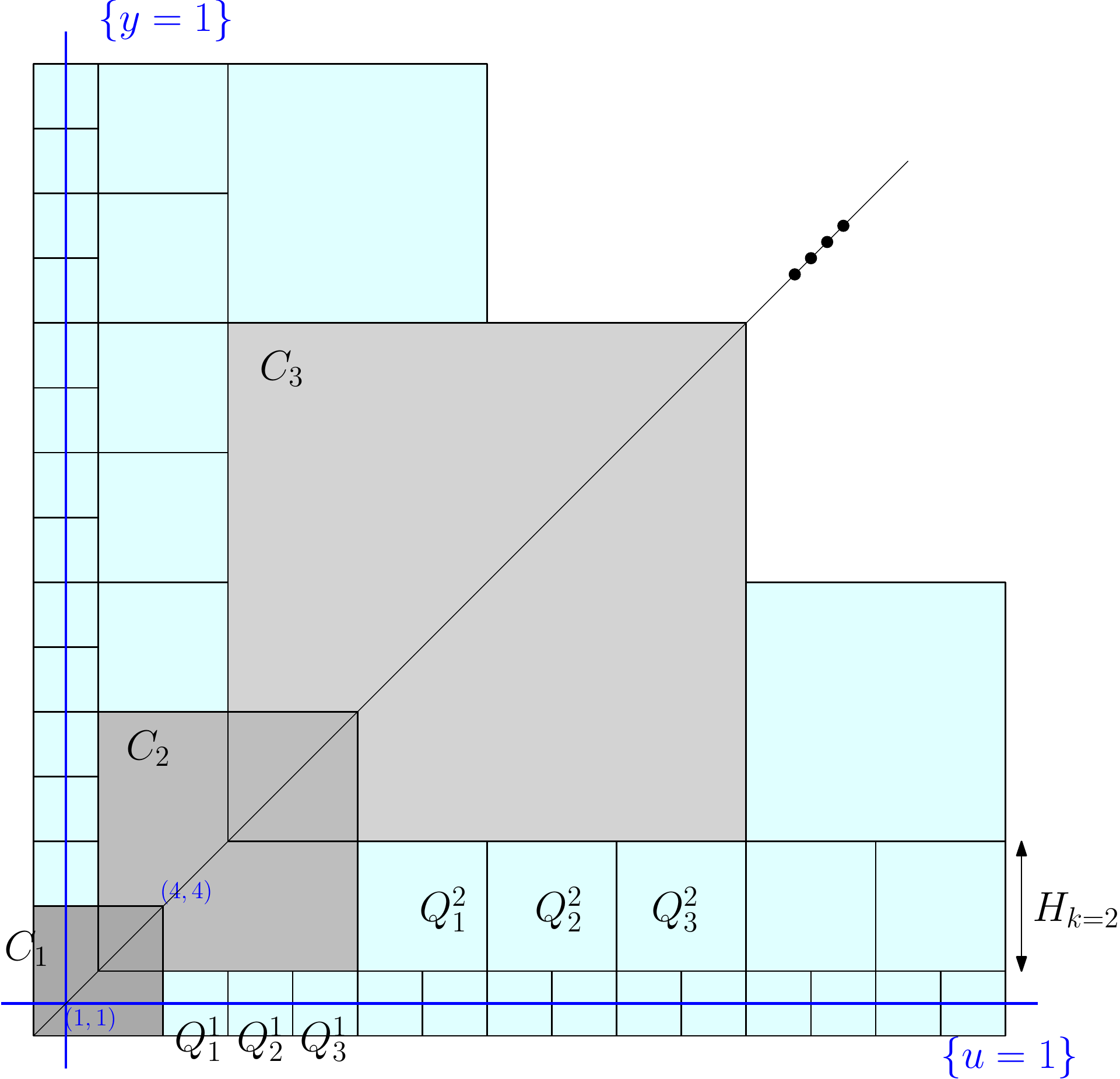}
\end{center}
\caption{We shall use the above decomposition into dyadic squares of $[0,\frac 2 \ell]^2$.
}\label{f.decomp}
\end{figure}

\ni
As a warm-up, let us analyse the easier case of $q=2$: this gives 

\begin{align*}\label{}
& \Eb{(\du_\eps)^2} \leq O(1) \Eb{ J_\epsilon } \\
& \leq O(1) 
\iint_{[-1,1]^2} 
\big(\frac{1}{|y+ \frac \ell 2|^{\tfrac{1}{2}-H}}
- \frac{1}{|y- \frac \ell 2|^{\tfrac{1}{2}-H}}
- \frac{1}{|u + \frac \ell 2|^{\tfrac{1}{2}-H}}
+ \frac{1}{|u - \frac \ell 2|^{\tfrac{1}{2}-H}}
\big)^2\frac{1}{|u-y|_\eps^{1-  2\tilde H  + 4 \gamma^2 }} du dy\,.
\end{align*}
We used here the fact that $\E[M_{2\gamma}^\epsilon(du)M_{2\gamma}^\epsilon(dy) ]=e^{4 \gamma^2 \widehat{C}_\epsilon(u-y)}dudy \leq O(1) 
|u-y|_\eps^{- 4 \gamma^2  }dudy$.
Then we change of scales as follows,
\begin{align*}\label{}
& \Eb{(\du_\eps)^2} \leq O(1)
\ell^{2(H+\tilde H) - 4 \gamma^2} \times\\
&\iint_{[-2/\ell,2/\ell]^2} 
\big(\frac{1}{|y +1|^{\tfrac{1}{2}-H}}
- \frac{1}{|y- 1|^{\tfrac{1}{2}-H}}
- \frac{1}{|u + 1|^{\tfrac{1}{2}-H}}
+ \frac{1}{|u - 1|^{\tfrac{1}{2}-H}}
\big)^2\frac{1}{|u-y|_{2\eps/\ell}^{1- 2\tilde H  +  4 \gamma^2}} du dy.
\end{align*}
We are thus left with studying the above integral. We shall focus on $[0,2/\ell]^2$ and rely on the decomposition of that square into dyadic squares defined as follows (See Figure \ref{f.decomp}):
\[
[0,2/\ell]^2 = \bigcup_{k\geq 1} C_k \cup H_k \cup V_k\,,
\]
where 
\bi
\item $C_1:=[0,4]^2$
\item More generally, for all $k\geq 1$, let 
$$C_k:=[2^k-2, 2^k-2+2^{k+1}]^2\,.$$
\item For all $k\geq 1$, let $H_k$ be the ``corridor'' on the right of $C_k$ as in Figure \ref{f.decomp}, i.e.
\[
H_k:= [2^k-2+2^{k+1}, \frac 2 \ell] \times [2^k-2, 2^k-2+ 2^{k}]
\]
And let $V_k$ the ``corridor'' on the top of $C_k$, i.e. 
\[
V_k:= [2^k-2, 2^k-2+2^{k}] \times [2^k-2+2^{k+1}, \frac 2 \ell] 
\]
\item We shall also split each corridor $H_k$ (and equivalently $V_k$ but by symmetry we will never analyze this case) into dydic squares $\{Q_m^k\}_{m =1,\ldots, \log_2 \frac {2^{-k}} {\ell} }$ of width $2^{k}$ as in figure \ref{f.decomp}. 
\item Finally, let us point out that this division is well adapted to the bottom/left corner of $[0,\frac 2 \ell]^2$ (which as we will see will give the main contributions to $\Eb{|\delta_\ell u|^q}$) but will not match nicely with the right and top boundaries of $[0,\frac 2 \ell]^2$. As the contributions of the squares $C_k$ and $Q_m^k$ will be shown to be negligible at that distance, we will not bother with adapting the shape of these limiting squares. 
\ei

\subsubsection{Analyzing the contribution of the square $C_1$}

There are singularities in this special square that need some care: $\{y=1\}$, $\{u=1\}$ and $\{u=y\}$. 


\begin{align*}\label{}
& 
\iint_{[0,4]^2} 
\big(\frac{1}{|y +1|^{\tfrac{1}{2}-H}}
- \frac{1}{|y- 1|^{\tfrac{1}{2}-H}}
- \frac{1}{|u + 1|^{\tfrac{1}{2}-H}}
+ \frac{1}{|u - 1|^{\tfrac{1}{2}-H}}
\big)^2\frac{1}{|u-y|_{2\eps/\ell}^{1- 2\tilde H + 4 \gamma^2  }} du dy  \\
& \leq 
C   \int_{[0,4]^2}
\big(
\frac{1}{|y- 1|^{\tfrac{1}{2}-H}}
- \frac{1}{|u - 1|^{\tfrac{1}{2}-H}}
\big)^2\frac{1}{|u-y|_{2\eps/\ell}^{1-  2\tilde H  +  4 \gamma^2  }} du dy  
\end{align*}

To analyze this integral, we will use a certain partition of $[0,4]^2$ and use scaling arguments around the triple singularity $x_0=(1,1)$.
Introduce for any $D\subset C_1=[0,4]^2$,
\[
A_\eps(D):= \int_{D}
\big(
\frac{1}{|y- 1|^{\tfrac{1}{2}-H}}
- \frac{1}{|u - 1|^{\tfrac{1}{2}-H}}
\big)^2\frac{1}{|u-y|_{2\eps/\ell}^{1- 2\tilde H  +  4 \gamma^2  }} du dy  
\]
Clearly, one has 
\[
A_\eps([0,4]^2) = A_\eps([0,2]^2) + A_\eps([0,4]^2\setminus [0,2]^2)
\]
The second term is easier to analyse and we shall focus on the first one which we decompose as follows:
\begin{align*}\label{}
A_\eps([0,2]^2) = & A_\eps([1/2,3/2]^2) + A_\eps([1/2,3/2]\times [0,2]\setminus [1/2,3/2]) \\
& + A_\eps( [0,2]\setminus [1/2,3/2] \times [1/2,3/2]) + \text{corners} .
\end{align*}
A straightforward scaling shows that $A_{\eps/2}([1/2,3/2]^2) = 2^{4\gamma^2 - 2(H+\tilde H)} A_\eps([0,2]^2)$. Furthermore, it is easy to see that the other terms listed above are uniformly bounded as $\eps\to 0$ and for the corners, the function to be integrated near the diagonal behaves like
\[
\big(\frac{1}{|y- 1|^{\tfrac{1}{2}-H}}
- \frac{1}{|u - 1|^{\tfrac{1}{2}-H}}
\big)^2\frac{1}{|u-y|_{2\eps/\ell}^{1- 2\tilde H  +  4 \gamma^2  }} 
\asymp (|u-y|)^2 \frac{1}{|u-y|_{2\eps/\ell}^{1- 2\tilde H  + 4 \gamma^2  }} 
\]
which gives the additional constraint $1+2\tilde H - 4\gamma^2 >-1$, i.e $2\gamma^2 < 1+\tilde H$, which is satisfied because $\tilde H \geq 0$ and $2\gamma^2<1$.

\subsubsection{Analyzing the contribution of the squares $\{C_k\}_{k\geq 2}$}

Let us introduce the following function

\begin{align*}\label{}
K(x):=  \left|\frac 1 {x+1}\right|^{\frac 1 2 -H} - \left|\frac 1 {x-1}\right|^{\frac 1 2 - H} \sim_{x\to \infty} c \left|\frac 1 {x}\right|^{3/2-H}
\end{align*}
It is straighforward to check that as $x,y \in C_k$, $k\geq 2$, one has 
\begin{align*}\label{}
|K(x) - K(y)| \leq  \|K'\|_{\infty, [2^k, 2^{k+1}]} |x-y| \leq C (2^{-k})^{\frac 5 2 - H} |x-y| 
\end{align*}

This implies that if we define 
\begin{align*}\label{}
h(y,u):= \big(\frac{1}{|y +1|^{\tfrac{1}{2}-H}}
- \frac{1}{|y- 1|^{\tfrac{1}{2}-H}}
- \frac{1}{|u + 1|^{\tfrac{1}{2}-H}}
+ \frac{1}{|u - 1|^{\tfrac{1}{2}-H}}
\big)^2\frac{1}{|u-y|_{2\eps/\ell}^{1-  2\tilde H  + 4 \gamma^2  }} 
\end{align*}
then one has 
\begin{align*}\label{}
\int_{C_k}h(y,u)\,dydu & \leq  C (2^{-k})^{5 -2H} \int_{C_k}  |u-y|^{1+2 \tilde H- 4\gamma^2}dydu \\
& \leq C (2^{-k})^{2- 2H - 2\tilde H + 4\gamma^2}.
\end{align*}
This handles the contribution given by the squares $\{C_k\}_{k\geq 2}$:
\begin{align*}\label{}
\ell^{2(H+\tilde H) - 4 \gamma^2}  \sum_{k= 2}^{\log_2 \frac 2 \ell} \iint_{C_k} h(y,u) dy du 
& \leq O(1) \ell^{2(H+\tilde H) - 4 \gamma^2}  \sum_{k\geq 1}  (2^{-k})^{2- 2H - 2\tilde H + 4\gamma^2}\\
& \leq O(1) \ell^{2(H+\tilde H) - 4 \gamma^2} 
\end{align*}
as we  $2H + 2 \tilde H - 4 \gamma^2 < 2$. 

\subsubsection{Analyzing the contribution of the corridors $\{H_k\}_{k\geq 1}$ and $\{V_k\}_{k\geq 1}$}

By symmetry, we will only focus on the horizontal corridors $\{H_k\}_{k\geq 1}$. Here $k$ ranges from 1 to $\log_2 \frac 1 \ell$.

The first corridor $k=1$ will need a separate study as it is traversed throughout by the line-singularity $\{u=1\}$.

\ni
\underline{Corridors $k=2, \ldots, \log_2 \frac 1 \ell$:} 
The horizontal corridor $H_k$ is made of the $2^k$-squares $Q_m^k$ with $m=1,
\ldots, \log_2 \frac{ 2^{-k}} {2 \ell}$. (See Figure \ref{f.decomp}). 

\begin{align*}
& \sum_{m=1}^{\log_2 \frac{ 2^{-k}} {2 \ell}}\iint_{Q_m^k} 
\big(\frac{1}{|y +1|^{\tfrac{1}{2}-H}}
- \frac{1}{|y- 1|^{\tfrac{1}{2}-H}}
- \frac{1}{|u + 1|^{\tfrac{1}{2}-H}}
+ \frac{1}{|u - 1|^{\tfrac{1}{2}-H}}
\big)^2\frac{1}{|u-y|_{2\eps/\ell}^{1-  2\tilde H  +  4 \gamma^2  }} du dy  \\
& \leq O(1) \sum_{m=1}^{\log_2 \frac{ 2^{-k}} {2 \ell}}
\mathrm{Area}(Q_m^k)  \| h \|_{\infty, Q_m^k} 
\end{align*}
Using the fact that $K(x)\sim_{x\to \infty} x^{-3/2+H}$, it is straightforward to check that 
\begin{align*}\label{}
 \| h \|_{\infty, Q_m^k} & \leq C \, [(2^{-k})^{3/2-H}]^2 \left(\frac 1 {m 2^k}\right)^{1- 2\tilde H + 4 \gamma^2} \\
 & \leq  C\,  [2^{-k}]^{4-2H-2\tilde H+4\gamma^2} m^{-1 + 2\tilde H - 4  \gamma^2},
\end{align*}
which gives us:
\begin{align*}\label{}
\sum_{m=1}^{\log_2 \frac{ 2^{-k}} {2 \ell}}
\mathrm{Area}(Q_m^k)  \| h \|_{\infty, Q_m^k} 
& \leq O(1) 2^{2k} [2^{-k}]^{4-2H-2\tilde H+4\gamma^2} \sum_{m=1}^{\log_2 \frac{ 2^{-k}} {2 \ell}} m^{-1 + 2\tilde H - 4  \gamma^2}  \\
& \leq O(1) [2^{-k}]^{2-2H-2\tilde H+4\gamma^2}  \sum_{m=1}^{\log_2 \frac{ 2^{-k}} {2 \ell}} m^{-1 + 2\tilde H - 4  \gamma^2}
\end{align*}
By  our assumptions above, we see that the exponent of $2^{-k}$ is negative. Our assumption $\tilde H < 2 \gamma^2$ entails that  we obtain an exponent $\alpha>1$ in $\sum_m m^{-\alpha}$.  This shows that  the main contribution will come from the first corridors.

\ni
\underline{First corridor ($k=1$)}

In this case the line-singularity $\{u=1\}$ traverses all squares $\{Q_m^{k=1}\}_{m=1,\ldots, \log_2 1/(4\ell)}$. This singularity is easier to deal with than the "2 lines singularity" for the above square $C_1$ as the singularity is integrable in $\int_0^2 du$. Therefore it is not hard to obtain the following upper bound on the first corridor $H_1$:
\begin{align*}\label{}
& \sum_{m=1}^{\log_2 \frac{ 1} {4 \ell}}\iint_{Q_m^1} 
\big(\frac{1}{|y +1|^{\tfrac{1}{2}-H}}
- \frac{1}{|y- 1|^{\tfrac{1}{2}-H}}
- \frac{1}{|u + 1|^{\tfrac{1}{2}-H}}
+ \frac{1}{|u - 1|^{\tfrac{1}{2}-H}}
\big)^2\frac{1}{|u-y|_{2\eps/\ell}^{1- 2\tilde H + 4 \gamma^2 }} du dy  \\
& \leq O(1)  \sum_{m=1}^{\log_2 \frac{ 1} {4 \ell}}\iint_{Q_m^1}  m^{-1 + 2\tilde H - 4\gamma^2} \leq O(1).
\end{align*}

Summarizing the above estimates, we thus obtain the following sharp bound (up to multiplicative constants)
\begin{align*}\label{}
\Eb{(\delta_\ell u)^2} \asymp  \ell^{2(H+\tilde H) - 4 \gamma^2} 
\end{align*}
which shows as expected that 
\begin{align*}\label{}
\xi(q=2) = 2(H+\tilde H) - 4 \gamma^2.
\end{align*}

\subsection{$q\geq 2$ moments and non-linear spectrum}\label{ss.NLS}
Recall the following estimate
\begin{align*}\label{}
& \Eb{|\du_\eps(x=0)|^q}  \leq \tilde C_q\Eb{J_\epsilon^{q/2}},
\end{align*}
where
$$ J_\epsilon=\iint_{[-1,1]^2} 
\big(\frac{1}{|y+ \frac \ell 2|^{\tfrac{1}{2}-H}}
- \frac{1}{|y- \frac \ell 2|^{\tfrac{1}{2}-H}}
- \frac{1}{|u + \frac \ell 2|^{\tfrac{1}{2}-H}}
+ \frac{1}{|u - \frac \ell 2|^{\tfrac{1}{2}-H}}
\big)^2\frac{1}{|u-y|_\eps^{1-  2\tilde H}}M_{2\gamma}^\epsilon(du)M_{2\gamma}^\epsilon(dy).$$
Let us focus on the square $[0,1]^2$ ($[-1,0]^2$ is treated the same way and the two other squares are easier to deal with).
By using the triangular inequality for the $L_{q/2}$-norm (we assume here that $q\geq 2$), we have the bound:
\begin{align*}\label{}
& \Eb{|\du_\eps(x=0)|^q}^{2/q}  \\
& \leq O(1)  \sum_{k=1}^{\log \frac 2 \ell}  
\sum_{A\in \{ \frac \ell 2 C_k, \frac \ell 2 V_k,\frac \ell 2 H_k\}} \left\| \iint_A  \ell^{-2 + 2H + 2\tilde H} g(\frac {2y} \ell  ,\frac {2u} \ell)
M_{2\gamma}^\epsilon(du)M_{2\gamma}^\epsilon(dy) \right\|_{q/2}
\end{align*}
where similarly as for the $q=2$ case, we shall use the function  
\begin{align*}\label{}
g(y,u):= \big(\frac{1}{|y +1|^{\tfrac{1}{2}-H}}
- \frac{1}{|y- 1|^{\tfrac{1}{2}-H}}
- \frac{1}{|u + 1|^{\tfrac{1}{2}-H}}
+ \frac{1}{|u - 1|^{\tfrac{1}{2}-H}}
\big)^2\frac{1}{|u-y|_{2\eps/\ell}^{1- 2\tilde H}} .
\end{align*}

Now, on $\R_+^2$, it is easy to check that 
\begin{align*}
g(y,u) & \leq  2 \Bigl[ (\frac{1}{|y +1|^{\tfrac{1}{2}-H}}
- \frac{1}{|u + 1|^{\tfrac{1}{2}-H}})^2 
+ 
(\frac{1}{|u - 1|^{\tfrac{1}{2}-H}}
- \frac{1}{|y- 1|^{\tfrac{1}{2}-H}}
)^2 
\Bigr]\frac{1}{|u-y|_{2\eps/\ell}^{1- 2\tilde H}}  \\
& \leq 4 (\frac{1}{|u - 1|^{\tfrac{1}{2}-H}}
- \frac{1}{|y- 1|^{\tfrac{1}{2}-H}}
)^2  \frac{1}{|u-y|_{2\eps/\ell}^{1- 2\tilde H}}  =:f_\eps(y,u)
\end{align*}

\subsubsection{Analyzing the contribution of the square $C_1$}

As in the easier case of $q=2$, we have to analyze the quantity
\begin{align*}\label{}
&  \left\| \iint_{\frac \ell 2 C_1} \ell^{-2 + 2H + 2\tilde H} f_\eps(\frac {2y} \ell  ,\frac {2u} \ell )
M_{2\gamma}^\epsilon(du)M_{2\gamma}^\epsilon(dy) \right\|_{q/2}
\end{align*}
We will proceed here as in Subsection \ref{ss.mom} by relying on Kesten's inequality as well as scaling arguments. 
The main singularities in the square $\frac \ell 2 C_1 = [0,2 \ell]^2$ arise in the sub-square $[0,\ell]^2$ on which we now focus. We will use the scaling properties of the chaos twice: first to take into account the highly correlated nature of $\hat X$ within $[0,\ell]^2$ and then a second time along a similar decomposition as in Subsection \ref{ss.mom} (except it will be centered here around the point $x_0=(\frac \ell 2, \frac \ell 2)$) to prove that the limit as $\eps \to 0$ exists. 

\medskip

\ni
\textbf{First rescaling.}
We need to estimate the quantity 
\begin{align*}\label{}
R_{C_1}&:=   \ell^{-2 + 2H + 2\tilde H} \Eb{\bigl(\iint_{[0,\ell]^2}  f_\eps(\frac {2y} \ell,\frac {2u} \ell)
M_{2\gamma}^\epsilon(du)M_{2\gamma}^\epsilon(dy)  \bigr)^{q/2}}^{2/q} \\
& =  \ell^{-2 + 2H + 2\tilde H} 
 \Eb{\bigl(\iint_{[0,\ell]^2}  f_\eps(\frac {2y} \ell,\frac {2u} \ell)
e^{2\gamma \widehat{X}_{\epsilon}(u) -2\gamma^2 c_{\epsilon}}
e^{2\gamma \widehat{X}_{\epsilon}(y) -2\gamma^2 c_{\epsilon}} du dy
\bigr)^{q/2}}^{2/q}.
\end{align*}

Using the change of variables $\bar y = y/\ell$, $\bar u = u/\ell$ together with the identity \eqref{exact}, we have 
\begin{align}\label{}
R_{C_1} & = 
 \ell^{2H + 2\tilde H} 
\Eb{\bigl(\iint_{[0,1]^2}  f_\eps(2 \bar y , 2 \bar u )
e^{2\gamma \widehat{X}_{\epsilon}(\ell \bar u) -2\gamma^2 c_{\epsilon}}
e^{2\gamma \widehat{X}_{\epsilon}(\ell \bar y) -2\gamma^2 c_{\epsilon}} du dy
\bigr)^{q/2}}^{2/q} \nn \\
& = 
 \ell^{2H + 2\tilde H}
\E\big[e^{2\gamma q\Omega_{1/\ell}}\big]^{2/q}e^{-4\gamma^2\ln 1/\ell}
\Eb{\bigl(\iint_{[0,1]^2}  f_\eps(2 \bar y , 2 \bar u )
M_{2\gamma}^{\epsilon/\ell}(d \bar u)M_{2\gamma}^{\epsilon/\ell}(d\bar y)
\bigr)^{q/2}}^{2/q} \nn \\
& = \ell^{2(H+\tilde H) + 4 \gamma^2 - 4 \gamma^2 q}
\Eb{\bigl(\iint_{[0,1]^2}  f_\eps(2 \bar y , 2 \bar u )
M_{2\gamma}^{\epsilon/\ell}(d \bar u)M_{2\gamma}^{\epsilon/\ell}(d\bar y)
\bigr)^{q/2}}^{2/q} .
\end{align}

\ni
\textbf{Second rescaling.}
Now using the same notations as in Subsection \ref{ss.mom},  define 
\begin{align*}\label{}
A_{\eps,\ell}(D):= \Eb{\bigl(\iint_D  f_\eps( 2\bar y ,2 \bar u)
M_{2\gamma}^{\epsilon/\ell}(du)M_{2\gamma}^{\epsilon/\ell}(dy)  \bigr)^{q/2}}^{2/q}
\end{align*}
Again by the triangle inequality for $L^p$ norms with $p=\frac q 2 \geq 1$, we have 

\begin{align*}\label{}
A_{\eps/2,\ell}([0, 1]^2) \leq  & A_{\eps/2,\ell}([\frac 1 4, \frac 3 4]^2) +  A_{\eps/2,\ell}([\frac 1 4, \frac 3 4]\times[0,1]\setminus [\frac 1 4, \frac 3 4]) \\
&  +  A_{\eps/2,\ell}([0,1]\setminus [\frac 1 4, \frac 3 4] \times [\frac 1 4, \frac 3 4]) + \text{the 4 corner squares} 
\end{align*}

Let us first deal with the most problematic square : the one centred around the point-singularity $(1/2,1/2)$.
\begin{lemma}\label{}
\begin{align*}\label{}
 A_{\eps/2,\ell}([\frac 1 4, \frac 3 4]^2)  \leq 2^{-2(H+\tilde H)-4\gamma^2+4\gamma^2 q}   A_{\eps,\ell}([0, 1]^2)
\end{align*}
\end{lemma}
\ni {\it Proof.} We will use the fact (already used above in \eqref{exact})
\begin{equation}\label{exactwithell}
(\widehat{X}_{\lambda\epsilon }(\lambda u + \lambda (\ell/2, \ell/2)))_{u\in [0,\ell]}=(\widehat{X}_{\epsilon}(u)+\Omega_\lambda)_{u\in [0,\ell]}
\end{equation}
where $\Omega_\lambda$ is a centered Gaussian random variable with variance $-\ln \lambda$   independent of the process $\widehat{X}_{\epsilon}$.
By making the change of variable $u'=2\bar u -1/2$ and $y'=2\bar y -1 /2$ we get 
\begin{align*}\label{}
A_{\eps/2,\ell}([\frac 1 4, \frac 3 4]^2)
& = \Eb{\bigl(\iint_{[\frac 1 4, \frac 3 4]^2}  f_{\eps/2}(2 \bar y, 2 \bar u)
e^{2\gamma \widehat{X}_{\epsilon/(2\ell)}(\bar u) -2\gamma^2 c_{\epsilon/(2\ell)}}
e^{2\gamma \widehat{X}_{\epsilon/(2\ell)}(\bar y) -2\gamma^2 c_{\epsilon/(2\ell)}} d\bar u d\bar y \bigr)^{q/2}}^{2/q} \\
& = \frac 1 4 
\E\big[e^{2\gamma q\Omega_{1/2}}\big]^{2/q}e^{-4\gamma^2\ln 2}  \\
& \times \Eb{\bigl(\iint_{[0, 1]^2}  f_{\eps/2}( 2 (y'/2 +1/4),2  (u'/2 +1/4))
M_{2\gamma}^{\epsilon/\ell}(du')M_{2\gamma}^{\epsilon/\ell}(dy') \bigr)^{q/2}}^{2/q} \\
& = 2^{-2(H+\tilde H)-4\gamma^2+4\gamma^2 q}  \,\, A_{\eps,\ell}([0, 1]^2)  
\end{align*}\qed

We thus need to assume that the exponent $4 \gamma^2 q - 4\gamma^2 -2 (H+\tilde H)$ is negative, or otherwise stated $q < 1+\frac 1 {2\gamma^2} (H+\tilde H)$.  
It thus remains to control the other squares in order to show the following estimate:
\begin{lemma}\label{l.4squares}
Assuming  $q < 1+\frac 1 {2\gamma^2} (H+\tilde H)$ and  $q<2(1+\frac H {\gamma^2})\wedge \frac 1 {2 \gamma^2}$
then, there exists constants $r<1$ and $C<\infty$ s.t. uniformly in $0<\eps<\ell$, 
\[
 A_{\eps/2,\ell}([0, 1]^2)  \leq r\, A_{\eps,\ell}([0, 1]^2) + C 
\]
\end{lemma}

\ni {\it Proof.} Let us focus for example on the rectangle $R_{top}:=[1/4,3/4]\times[3/4,1]$. (The diagonal squares  are less singular and are treated as in Subsection \ref{ss.mom}). The exact same analysis as the one carried for the square $[0,1/2] \times [1/2,1]$ in Subsection \ref{ss.mom} applies here given the additional constraint that $q<2(1+\frac H {\gamma^2})\wedge \frac 1 {2 \gamma^2}$.  \qed

All together, we see that under the conditions of Lemma \ref{l.4squares}, the contribution of the square $\frac \ell 2 C_1$ to the $q^{th}$ moment of the increment (to the power $2/q$) is given by 
\[
O(1) \ell^{2(H+\tilde H) + 4 \gamma^2 - 4\gamma^2 q}
\]
which is what we wanted.

\subsubsection{Analyzing the contribution of the squares $\{C_k\}_{k\geq 2}$}

Easiest case: here no singularities, and standard scaling argument.

\subsubsection{Analyzing the contribution of the corridors $\{H_k\}_{k\geq 1}$ and $\{V_k\}_{k\geq 1}$}

Let us start with the following slight generalisation of the celebrated Kahane's convexity inequality. This small extension is of independent interest. 

\begin{proposition}\label{pr.KahExt}
Let $(X_i)_{1\leq i \leq n+m}$ and $(Y_i)_{1\leq i \leq n+m}$ be two centered Gaussian vectors satisfying for all $i,j$
\begin{align*}\label{}
\Eb{X_i X_j} \leq \Eb{Y_i Y_j}\,.
\end{align*}
Then for all sequence of nonnegative weights $(p_i)_{1\leq i \leq n+m}$ and all \textbf{increasing} convex functions $F: \R_+ \to \R_+$ and $G : \R_+ \to \R_+$, one has 
\begin{align*}\label{}
\Eb{F(\sum_{i=1}^n p_i e^{X_i - \frac 1 2 \Eb{X_i^2}} ) G(\sum_{k=n+1}^m p_k e^{X_k - \frac 1 2 \Eb{X_k^2}} )}
& \leq 
\Eb{F(\sum_{i=1}^n p_i e^{Y_i - \frac 1 2 \Eb{Y_i^2}} ) G(\sum_{k=n+1}^m p_k e^{Y_k - \frac 1 2 \Eb{Y_k^2}} )}
\end{align*}
\end{proposition}

\ni {\it Proof.}  The proof follows exactly the same lines as the original proof by Kahane (see for example \cite{Kah85, RobVar08}). We shall only sketch briefly how to adapt the proof here. Consider two independent realizations of the Gaussian vectors 
$(X_i)$ and $(Y_i)$ and interpolate between the two as follows:
\[
Z_i(t):= \sqrt{t} X_i + \sqrt{1-t} Y_i.
\]
Consider the function 
\begin{align*}\label{}
\phi(t):= \Eb{F(\sum_{i=1}^n p_i e^{Z_i(t) - \frac 1 2 \Eb{Z_i(t)^2}} ) G(\sum_{k=n+1}^m p_k e^{Z_k(t) - \frac 1 2 \Eb{Z_k(t)^2}} )}.
\end{align*}

Then, by using the Gaussian integration by parts formula, it is not difficult to obtain  the following identity 
\begin{align*}\label{}
\phi'(t) 
& = \sum_{i,j=1}^n p_i p_j \left(( \Eb{X_iX_j} - \Eb{Y_i Y_j} )\Eb{e^{Z_i(t) + Z_j(t) - \frac 1 2 \Eb{Z_i(t)^2} - \frac 1 2 \Eb{Z_j(t)^2} } F''(V_{n,t})G(W_{m,t})}\right) \\
& + 2 \sum_{\substack{i=1,\ldots,n \\ k=n+1,\ldots,n+m}}
p_i p_k \left( (\Eb{X_iX_k} - \Eb{Y_i Y_k}) \Eb{e^{Z_i(t) + Z_k(t) - \frac 1 2 \Eb{Z_i(t)^2} - \frac 1 2 \Eb{Z_k(t)^2} } F'(V_{n,t})G'(W_{m,t})}\right) \\
& + \sum_{k,l=n+1}^{n+m} p_k p_l \left(( \Eb{X_kX_l} - \Eb{Y_k Y_l}) \Eb{e^{Z_k(t) + Z_l(t) - \frac 1 2 \Eb{Z_k(t)^2} - \frac 1 2 \Eb{Z_l(t)^2} } F(V_{n,t})G''(W_{m,t})}\right) 
\end{align*}
where 
\[
\begin{cases}
V_{n,t}& := \sum_{i=1}^n p_i e^{Z_i(t) - \frac 1 2 \Eb{Z_i(t)^2}} \\
W_{m,t}& := \sum_{k=n+1}^{n+m} p_k e^{Z_k(t) - \frac 1 2 \Eb{Z_k(t)^2}} 
\end{cases}
\]
With our above assumptions, it implies that $\phi'(t) \leq 0$ $ \forall t\in [0,1]$, which concludes our proof as this shows that $\phi(0)\geq \phi(1) $.  \qed 

\medskip
As in Kahane's work (\cite{Kah85}), the same inequality for continuous multiplicative chaos measures immediately follows.
\medskip

We will need the following decorrelation Lemma.
\begin{lemma}\label{l.decor}
Let $A,B$ be two disjoint intervals of length $|A|=|B|=u$ on $[-2,2] \subset \R$ s.t. that $\mathrm{dist}(A,B) \geq u$. Then we have for all $q<q^*= 2 \frac 2 {(2 \gamma)^2} = \frac 1 {\gamma^2}$, 
\begin{align*}\label{}
\Eb{M_{2\gamma}(A)^{q/2} M_{2\gamma}(B)^{q/2}} & \leq 
 \leq C   (\frac 1 {\mathrm{dist}(A,B)})^{ \gamma^2 q^2}\;\; 
 |A|^{(1+ 2 \gamma^2) q - \gamma^2  q^2} 
\end{align*}
\end{lemma}

\ni
\ni {\it Proof.} Recall that our log-correlated stationary centered Gaussian field, $\widehat{X}$, has the following covariance structure for all $\eps>0$ (where $\widehat{X}_\eps:= \rho_\eps* \widehat{X}$):
\begin{align}\label{e.COV2}
\widehat{C}_\epsilon(x)= \E[\widehat{X}_\epsilon(x)\widehat{X}_\epsilon(0)]\sim  \ln_+\frac{1}{|x|+\epsilon}
\end{align}
Let us construct the following centered Gaussian Field $Y(x)$ on $\R$: define
\begin{align*}\label{}
Y(x) := Z(x) + \lambda \mathcal{N}(0,1)
\end{align*}
where $Z(x)$ is the centered Gaussian log-correlated field on $\R$ with covariance kernel 
\[
\Cov(Z(x), Z(y)) = \log_+ \frac {\mathrm{dist}(A,B)} { 2 |x-y|}
\]
and $\mathcal{N}(0,1)$ is a global Gaussian variable independent of $Z$. See for example \cite{RhoVar14} for a discussion on the $\log_+$ covariance kernel. 
 Note that by construction, the field $Z(x)$ is independent of $Z(y)$ as soon as $|x-y| > \mathrm{dist}(A,B)/2$. Therefore $Z_{|A}$ and $Z_{|B}$ are independent log-correlated fields.

Let us now fix $\lambda$ so that 
\[
\lambda^2= \log_+ \frac 1 {\mathrm{dist}(A,B)} + K
\]
where $K$ is a large constant to be fixed later which will not depend on $A$ nor $B$. 
By the covariance structure~\eqref{e.COV2}, it readily follows that for all $x,y \in A\cup B$, 
\begin{align*}\label{}
\Eb{\widehat{X}_\eps(x) \widehat{Y}_\eps(y)} 
& \leq \Eb{\widehat{X}(x) \widehat{X}(y)} \\
& \leq \log_+ \frac 1 {|x-y|} + C_1 \\
& \leq \log_+ \frac {\mathrm{dist}(A,B)} { 2 |x-y|} + \log_+ \frac 1 {\mathrm{dist}(A,B)}  + C_2\\
&= \Cov(Z(x), Z(y)) + \lambda^2 = \Cov(Y(x), Y(y))\,.
\end{align*}
(We thus choose $K$ to be the constant $C_2$ in the third line). 

We are now in position to apply Proposition \ref{pr.KahExt} (or rather its straightforward extension to continuous multiplicative chaos). It gives 
\begin{align*}\label{}
& \limsup_{\eps \to 0} 
\Eb{M_{2\gamma}^\eps(A)^{q/2} M_{2\gamma}^\eps(B)^{q/2}} \\
& \leq \Eb{M_{2\gamma}(A)^{q/2} M_{2\gamma}(B)^{q/2}} \\
& \leq \Eb{e^{2\gamma q \mathcal{N}(0,\lambda^2) - \frac{(2\gamma)^2} 2 q \lambda^2 }}
\Eb{ M_{2\gamma}^Z(A)^{q/2}} \Eb{M_{2\gamma}^Z(B)^{q/2}} \\
& \leq O_K(1)  (\frac 1 {\mathrm{dist}(A,B)})^{2 \gamma^2 q^2 - 2 \gamma^2 q}\;\;\Eb{ M_{2\gamma}^Z(A)^{q/2}} \Eb{M_{2\gamma}^Z(B)^{q/2}}\,, 
\end{align*}
where $M_{2\gamma}^Z$ stands for the multiplicative chaos measure of exponent $2\gamma$ induced by the field $Z$. As this field $Z$ is $\log_+$ correlated below scales of width $\mathrm{dist}(A,B)$ it enjoys exact scaling relations. (See \cite{RhoVar14}).  Now standard scaling arguments for these measures (recall we assumed $u=|A| \leq \mathrm{dist}(A,B)$) give us the bound below at least if $q$ is not too large, namely $q/2 < \frac 2 {(2\gamma)^2}$ i.e. $q<q^*$ as stated in the Lemma. 
\begin{align*}\label{}
\Eb{ M_{2\gamma}^Z(A)^{q/2}} 
& \asymp  \Eb{e^{2\gamma \frac q 2 \mathcal{N}(0,\log \frac {\mathrm{dist}(A,B)} {|A|}) - \frac{(2\gamma)^2} 2 \frac q 2 \log \frac {\mathrm{dist}(A,B)} {|A|} }} |A|^{q/2}\\
& \asymp \left( \frac {\mathrm{dist}(A,B)} {|A|} \right)^{\frac {\gamma^2} 2 q^2 -\gamma^2 q} |A|^{q/2}.
\end{align*}

Combining the above two estimates, we find as expected
\begin{align*}\label{}
&  \Eb{M_{2\gamma}(A)^{q/2} M_{2\gamma}(B)^{q/2}}   \\
& \leq O_K(1) (\frac 1 {\mathrm{dist}(A,B)})^{2 \gamma^2 q^2 - 2 \gamma^2 q}\;\; 
\left( \frac {\mathrm{dist}(A,B)} {|A|} \right)^{\gamma^2  q^2 - 2 \gamma^2 q} |A|^q \\
& \leq O_K(1)   (\frac 1 {\mathrm{dist}(A,B)})^{ \gamma^2 q^2}\;\; 
 |A|^{(1+ 2 \gamma^2) q - \gamma^2  q^2} 
\end{align*}

\qed

\ni
\underline{Corridors $k=2, \ldots, \log_2 \frac 1 \ell$}

The horizontal corridor $\ell H_k$  (we do not zoom here by a factor of $\frac 1 
\ell$) is made of the $2^k$-squares $\ell Q_m^k$ with $m=1,
\ldots, \log_2 \frac{ 2^{-k}} {2 \ell}$. 

\begin{align*}\label{}
& \left\| \iint_{\ell H_k}  \ell^{-2 + 2H + 2\tilde H} g(\frac y \ell,\frac u \ell)
M_{2\gamma}^\epsilon(du)M_{2\gamma}^\epsilon(dy) \right\|_{q/2} \\
& \leq   \ell^{-2 + 2H + 2\tilde H} \sum_{m=1}^{\log_2 \frac{ 2^{-k}} {2 \ell}}
 \left\| 
\iint_{\ell Q_m^k} 
 g(\frac y \ell,\frac u \ell)
M_{2\gamma}^\epsilon(du)M_{2\gamma}^\epsilon(dy) 
 \right\|_{q/2}\\
& \leq O(1) \ell^{-2 + 2H + 2\tilde H} 
\sum_{m=1}^{\log_2 \frac{ 2^{-k}} {2 \ell}}
 \| g \|_{\infty, Q_m^k}  \Eb{M_{2\gamma}(I^k_0)^{q/2} M_{2\gamma}(I^k_m)^{q/2}}^{2/q}\,,
\end{align*}
where the intervals $I^k_0$ and $\{ I^k_m\}_{m\geq 1}$ are such that each dyadic square $\ell Q_m^k =  I^k_m \times  I^k_0$. These intervals are of length $| I^k_0| =| I^k_m| = \ell  2^k$ and are at distance $\mathrm{dist}( I^k_0,  I^k_m)=m \ell 2^k$ from each other. Therefore we obtain from Lemma \ref{l.decor} that 
\begin{align*}\label{}
& \left\| \iint_{\ell H_k}  \ell^{-2 + 2H + 2\tilde H} g(\frac y \ell,\frac u \ell)
M_{2\gamma}^\epsilon(du)M_{2\gamma}^\epsilon(dy) \right\|_{q/2} \\
& \leq O(1) \ell^{-2 + 2H + 2\tilde H} 
\sum_{m=1}^{\log_2 \frac{ 2^{-k}} {2 \ell}}
 \| g \|_{\infty, Q_m^k}  \Eb{M_{2\gamma}(I^k_0)^{q/2} M_{2\gamma}(I^k_m)^{q/2}}^{2/q} \\
 & \leq O(1) \ell^{-2 + 2H + 2\tilde H} 
\sum_{m=1}^{\log_2 \frac{ 2^{-k}} {2 \ell}}
 \| g \|_{\infty, Q_m^k} \left[   (\frac 1 {\mathrm{dist}(I^k_0,I^k_m)})^{ \gamma^2 q^2}\;\; 
 |I^k_0|^{(1+ 2 \gamma^2) q - \gamma^2  q^2}  \right]^{2/q} \\
&  \leq O(1) \ell^{-2 + 2H + 2\tilde H} 
\sum_{m=1}^{\log_2 \frac{ 2^{-k}} {2 \ell}}
 \| g \|_{\infty, Q_m^k} \frac 1 {m^{2 \gamma^2 q}} \left[   (\ell 2^k)^{- \gamma^2 q^2} (\ell 2^k)^{(1+ 2 \gamma^2) q - \gamma^2  q^2}  \right]^{2/q}  \\
 & \leq 
 \ell^{-2 + 2H + 2\tilde H} 
\sum_{m=1}^{\log_2 \frac{ 2^{-k}} {2 \ell}}
 \| g \|_{\infty, Q_m^k} \frac 1 {m^{2 \gamma^2 q}}  (\ell 2^k)^{(1+2 \gamma^2) 2 - 4\gamma^2 q }\end{align*}
 
Now, as in the case of $q=2$ (where we relied on the function $h$ rather than $g$) and using the fact that $\phi(x)\sim_{x\to \infty} x^{-3/2+H}$, it is straightforward to check that 
\begin{align*}\label{}
 \| g \|_{\infty, Q_m^k} & \leq C \, [(2^{-k})^{3/2-H}]^2 \left(\frac 1 {m 2^k}\right)^{1- 2\tilde H} \\
 & \leq  C\,  [2^{-k}]^{4-2H-2\tilde H} m^{-1 + 2\tilde H},
\end{align*}
which gives us:
\begin{align*}\label{}
& \left\| \iint_{\ell H_k}  \ell^{-2 + 2H + 2\tilde H} g(\frac y \ell,\frac u \ell)
M_{2\gamma}^\epsilon(du)M_{2\gamma}^\epsilon(dy) \right\|_{q/2} \\
& 
\leq O(1) \ell^{-2 + 2H + 2\tilde H}  (\ell 2^k)^{(1+2 \gamma^2) 2 - 4\gamma^2 q } 
 [2^{-k}]^{4-2H-2\tilde H}
\sum_{m=1}^{\log_2 \frac{ 2^{-k}} {2 \ell}}
\frac 1 {m^{1- 2 \tilde H + 2 \gamma^2 q}}  \\
& \leq O(1) \ell^{2H + 2 \tilde H + 4\gamma^2 - 4 \gamma^2 q}  [2^{-k}]^{2-2H-2\tilde H-4\gamma^2+4\gamma^2 q}\,,
\end{align*}
where by our assumption 3 above ($\tilde H < q\gamma^2$), this is indeed summable in $m$ for all exponent $q\geq 2$.  
Also, by our assumption 2 above ($2 H + 2\tilde H - 4\gamma^2 <2$), we see that the exponent of $2^{-k}$ is indeed positive for all exponents $q\geq 2$.  This shows that  the main contribution to the $\| \cdot \|_{q/2}$ norm will come from the first corridors and is given  by 
\[
O(1)  \ell^{2H + 2 \tilde H + 4\gamma^2 - 4 \gamma^2 q} .
\]

As such the contribution of the corridors $\ell H_k, k\geq 2$ to $\Eb{\| \delta_\ell u\|^q}$ will be of of order (after taking to the exponent $q/2$ above)
\[
O(1)  \ell^{qH + q\tilde H + 2q \gamma^2 - 2 q^2 \gamma^2 } 
\]
which is indeed our expected structure function $\xi(q) = q(H +\tilde H) + 2q \gamma^2 - 2 q^2 \gamma^2$.

\medskip

\ni
\underline{First corridor ($k=1$)}
Finally, it can be shown that under the same set of constraints, this corridor contributes also  $O(1)\ell^{2H + 2 \tilde H + 4\gamma^2 - 4 \gamma^2 q}$ to $\Eb{\| \delta_\ell u\|^q}^{2/q}$. We leave the details to the reader as this case in some sense interpolates between the square $C_1$ (three-lines singularities) and the corridors $\{H_k\}_{k\geq 2}$ (no line singularities but large width): indeed $H_1$ has one line-singularity throughout (see Figure \ref{f.decomp}) and large width. \qed

\appendix

\section{Numerics}

\subsection{Estimation of the function $f_H(h)$ and its sign}\label{ann:NumEstF_H}

\begin{figure}[h]
\begin{center}
\includegraphics[width=12cm]{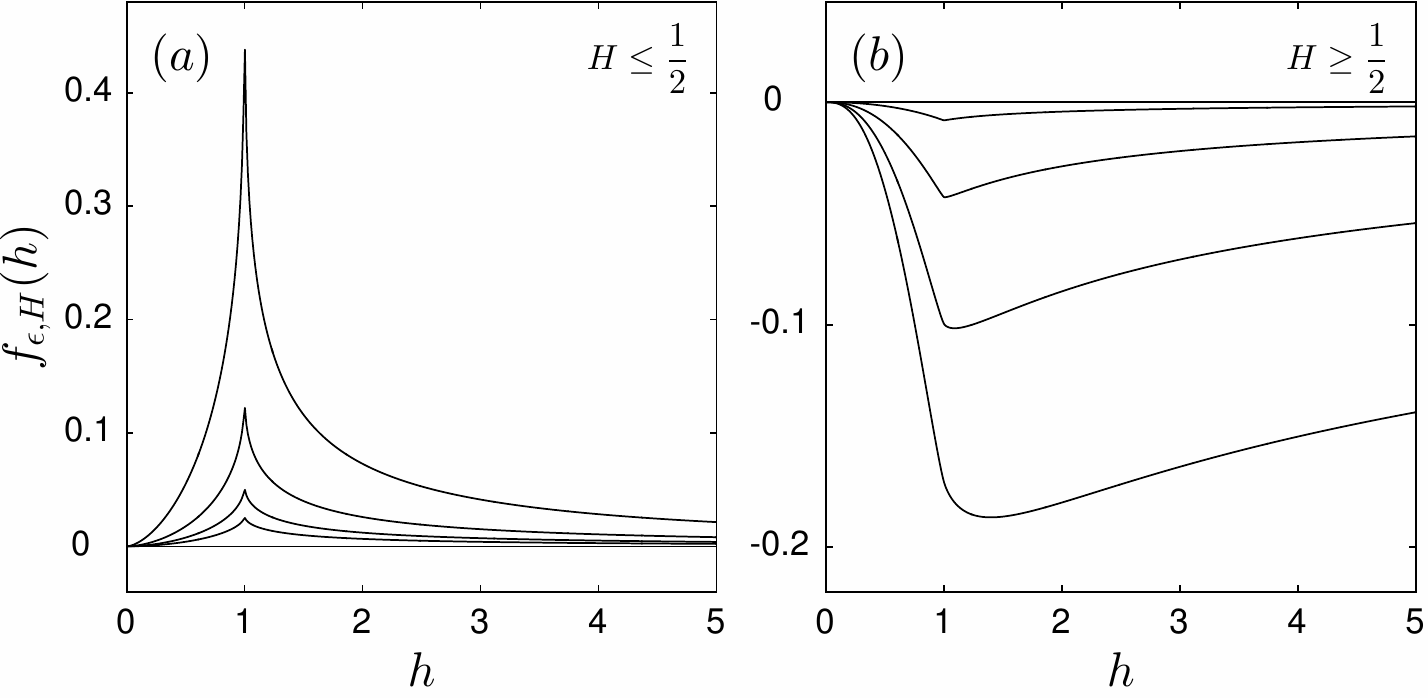}
\end{center}
\caption{\label{fig:f_H} Numerical estimation of the function $f_{\epsilon,H}(h)$, at a given approximation $\epsilon=10^{-5}$ (see text). (a) From top to bottom, $H=0.3, 0.35, 0.38, 0.4, 0.5$. (b) From top to bottom, $H=0.5, 0.6, 0.7,  0.8, 0.9$.}
\end{figure}

We represent in figure \ref{fig:f_H} the results of the numerical integration of a approximation $f_{\epsilon,H}$  of the function $f_{H}$ (Eq. \ref{deff_H}) entering in the third moment of the increments, namely
\begin{align}\label{deff_epsilonH}
 f_{\epsilon,H}(h)=\int &\left[\frac{1}{|x+1/2|_\epsilon^{\frac{1}{2}-H}}-\frac{1}{|x-1/2|_\epsilon^{\frac{1}{2}-H}}\right] \notag \\
 &\times\left[\frac{1}{|x+h+1/2|_\epsilon^{\frac{1}{2}-H}}-\frac{1}{|x+h-1/2|_\epsilon^{\frac{1}{2}-H}}\right]^2dx,
\end{align}
where enters the regularized norm $|x|_\epsilon^2 = |x|^2+\epsilon^2$. The numerical integration is made using adaptively a Newton-Cotes 5/9 point rule, as described in Ref. \cite{BerEsp91}. For this estimation, we use $\epsilon=10^{-5}$, and we checked (data not shown) that it is representative of the limit value $\epsilon\to 0$. We study a large set of values for $H$, and check that indeed $f_{\epsilon,1/2}(h)=0$ for $h>0$. This numerical study confirms the assumption on the sign of $f_H$ made in Eq. \ref{conjf_H}, that is positive for $H<1/2$, and negative for $H>1/2$. 

\subsection{Simulation of the random process, and estimation of its statistical properties}\label{ann:RandomSim}

\begin{figure}[h]
\begin{center}
\includegraphics[width=12cm]{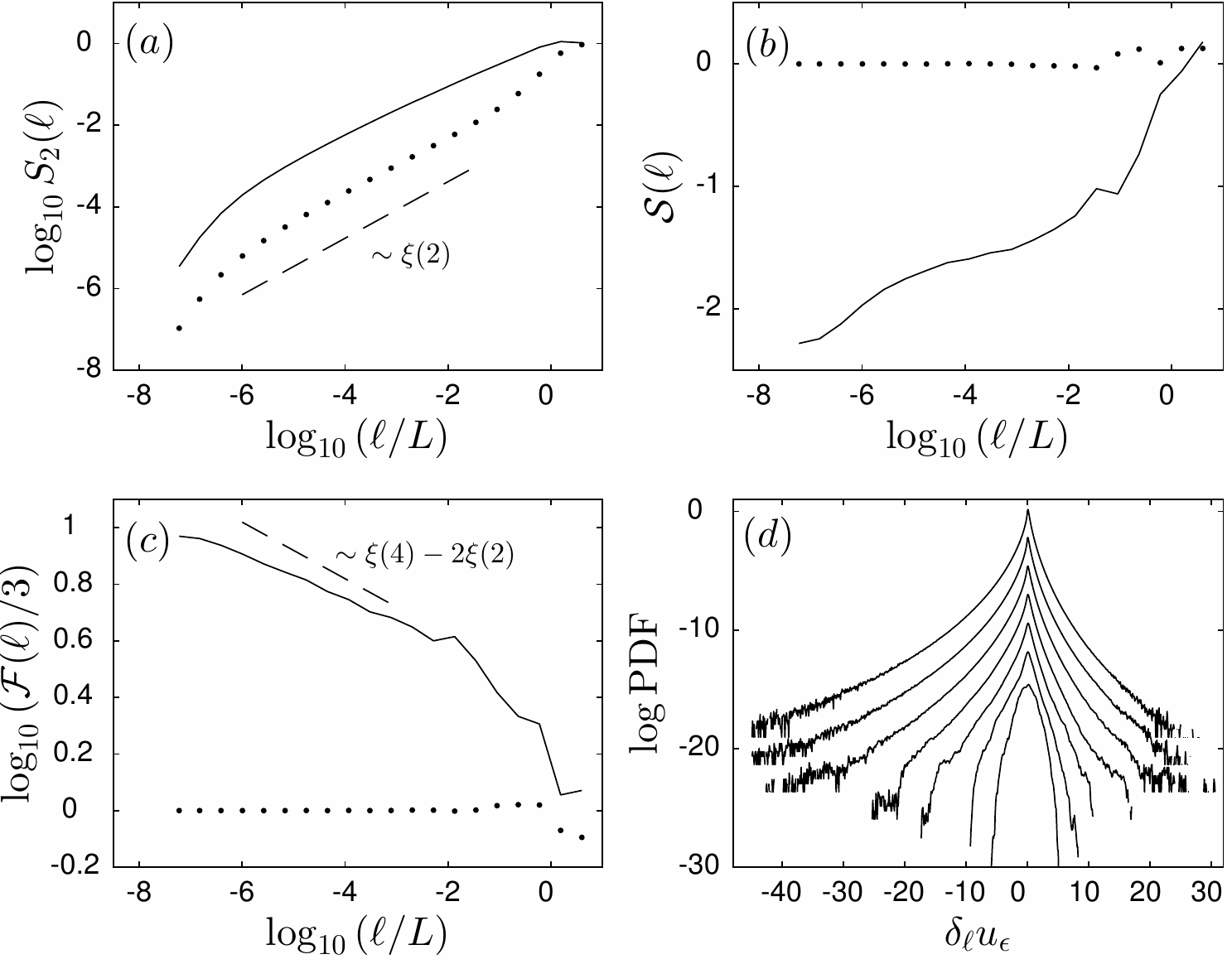}
\end{center}
\caption{\label{fig:StatProc} Numerical estimation of the statistical properties of $u_\epsilon$ \eqref{eq:SkMulProc} (using continuous lines), numerical details are provided in the text. (a) Estimation of  $S_2(\ell) = \E \left(\delta_\ell u_\epsilon \right)^2$ as a function of the scale $\ell$, in a logarithmic fashion. (b) Estimation of the skewness $\mathcal S(\ell)$ (Eq. \ref{eq:DefSkewnessNum}). (c) Estimation of the Flatness $\mathcal F(\ell)$ (Eq. \ref{eq:DefFlatNum}). In (a), (b) and (c), we superimpose the estimations of these statistical quantities using $u_\epsilon^{\mbox{\tiny g}}$ defined in \eqref{eq:ugepsilon} (dotted line) instead of $u_\epsilon$ (represented with a continuous line). In (a) and (c), we furthermore represent the expected power-law behaviors using dashed lines. (d) Logarithmic representations of the probability density functions (PDFs) of the increments $\delta_\ell u_\epsilon$ (renormalized by their respective standard deviation). Curves are arbitrary shifted vertically for clarity. From top to bottom, we have used $\log_{10}(\ell/L)=-7.2, -6, -4.8, -3.5, -2.3, -1,  0.2$. }
\end{figure}

We here present a method to simulate the proposed random field $u_\epsilon$ defined in \eqref{eq:SkMulProc} in a periodic fashion, such that we can work with the discrete Fourier transform. To do so, discretize the interval $[0,1]$ over $N$ collocation points. For full benefit of the fast Fourier transform (FFT) algorithm, choose $N$ to be a power of 2. This defines the numerical resolution of the simulation, i.e. $dx=1/N$. Choose for example as a cut-off function $\varphi_L $ a Gaussian shape, i.e. $\varphi_L(x) = e^{-\frac{x^2}{2L^2}}$. The precise shape of this function only matters at large scales, statistics at small scales are independent on it, besides its value at the origin. Choose as a regularized norm $|x|_\epsilon^2 = |x|^2+\epsilon^2$. Once again, the precise definition of the regularized norm does not matter since Theorem \ref{th1secdef} ensures that the statistical properties of $u_\epsilon$ are independent of the regularization procedure when $\epsilon\to 0$. 

Consider then two independent white fields $W$ and $\widehat{W}$, each of them made of $N$ independent realizations of a zero-average Gaussian variable of variance $dx$. Define the deterministic kernels $\phi_\epsilon$ \eqref{eq:DefphiIntro} (replace the norm $|.|$ entering in $\phi$ by its regularized form $|.|_\epsilon$) and $k_\epsilon$ \eqref{eq:DefkepsilonIntro} in a periodic fashion. Take $\widehat{X} = k_\epsilon\ast \widehat{W}$, and use $W$ as the remaining white field entering in the construction of $u_\epsilon$ \eqref{eq:SkMulProc}. Convolutions are then efficiently performed in the Fourier space.

We represent in Fig. \ref{fig:ShowProc} an instance of the process  $u_\epsilon$ \eqref{eq:SkMulProc}, as obtained by the aforementioned numerical method. We have used for the simulation the following set of parameters: $N= 2^{20}$, $L=1/3$, $\epsilon=2dx$, and the values $\gamma=\sqrt{0.025}/2$ and $H=1/3+4\gamma^2$. These chosen values for the parameters $\gamma$ and $H$ correspond to what is observed in turbulence (see Remark \ref{rem:RemTurbuIntro}).

To go further in the characterization of the statistical properties of the field $u_\epsilon$, we perform an additional simulation at a higher resolution $N=2^{31}$ in order to estimate in a reliable way the behaviors at small scales, and represent in Fig. \ref{fig:StatProc}  the results of our estimations. We have chosen  as a cut-off length scale $L=2^{-6}$, and as a regularizing scale $\epsilon=2dx$. Once again, values of the parameters are those which are realistic of turbulence, i.e. $\gamma=\sqrt{0.025}/2$ and $H=1/3+4\gamma^2$ (see Remark \ref{rem:RemTurbuIntro}). For the sake of comparison, we have furthermore made our estimations on the underlying fractional Gaussian field $u_\epsilon^{\mbox{\tiny g}}$ defined in \eqref{eq:ugepsilon}.


We begin in Fig. \ref{fig:StatProc}(a) with the estimation of the second-order structure function $S_2(\ell) = \E \left(\delta_\ell u_\epsilon \right)^2$ as a function of the scale. We observe for $u_\epsilon$ a power-law behavior at small scales, i.e. $S_2(\ell)\sim \ell^{\xi(2)}$, where the function $\xi(q)$ is defined in \eqref{xi}, consistently with Theorem \ref{th1secdef}. As for $u_\epsilon^{\mbox{\tiny g}}$, we also observe a power-law behavior at small scales, that we know to be $S_2(\ell)\sim \ell^{2H}$. Let us remark that since $\xi(2)$ is very close to $2H$, it is difficult to see a difference in between these power-laws. 

We represent in Fig. \ref{fig:StatProc}(b) the result of our estimation for the Skewness factor $\mathcal S(\ell)$ of the increments given by 
\begin{equation}\label{eq:DefSkewnessNum}
\mathcal S(\ell) = \frac{\E \left(\delta_\ell u_\epsilon \right)^3}{\left[ \E \left(\delta_\ell u_\epsilon \right)^2\right]^{3/2}}.
\end{equation}
We see that the present process is indeed skewed at small scales, being close to zero close to the large scale $L$, and growing towards values close to -2 at small scales. Remark that the quantity $\mathcal S(\ell)$is expected to behave as a power-law of exponent $\xi(3)-\frac{3}{2}\xi(2)$ at small scales. Remark also that it is indeed negative, as required by the phenomenology of turbulence (Section \ref{Physicsmotivation}). In comparison, we see that the Skewness factor for the Gaussian process $u_\epsilon^{\mbox{\tiny g}}$ is close to zero at any scales, as expected from symmetric statistical laws.

We represent in Fig. \ref{fig:StatProc}(c) the result of our estimation for the Flatness factor $\mathcal F(\ell)$ of the increments given by 
\begin{equation}\label{eq:DefFlatNum}
\mathcal F(\ell) = \frac{\E \left(\delta_\ell u_\epsilon \right)^4}{\left[ \E \left(\delta_\ell u_\epsilon \right)^2\right]^{2}}.
\end{equation}
Whereas $\mathcal F(\ell)$ is independent on the scale for the Gaussian process $u_\epsilon^{\mbox{\tiny g}}$ (and equal to 3), we see that it behaves as a power-law of exponent $\xi(4)-2\xi(2)$ at small scales.

Finally, we represent in Fig. \ref{fig:StatProc}(d) the histograms of the values of the increments $\delta_\ell u_\epsilon$ for several scales given in the caption. We see that whereas the histogram of the increments of the process $u_\epsilon$ are close to a Gaussian function at large scales $\ell \sim L$, they develop heavier and heavier tails at smaller scales, with a noticeable asymmetry.

\bibliographystyle{siam}


\end{document}